\documentclass[reqno,11pt,english]{amsart}
\usepackage[T1]{fontenc}
\usepackage[latin9]{inputenc}
\usepackage{color}
\usepackage{amstext}
\usepackage{amsthm}
\usepackage{amssymb}
\usepackage{comment}
\makeatletter
\usepackage{amsmath,amsfonts,amssymb,amsthm,epsfig}

\usepackage{color}
\usepackage[final,colorlinks,linkcolor=blue,anchorcolor=red,citecolor=blue]{hyperref}

\def\R{\mathbb R}
\def\N{\mathbb N}

\def\al{\alpha}
\def\be{\beta}
\def\ga{\gamma}
\def\de{\delta}
\def\ep{\epsilon}
\def\la{\lambda}

\def\ta{\theta}
\def\var{\varphi}
\def\om{\omega}
\def\na{\nabla}
\def\Om{\Omega}  
\def\De{\Delta}      

\def\cal{\mathcal}
\def\wq{\infty}
\def\pa{\partial}
\def\divergence{\text{\rm div}\,}

\def\loc{\text{\rm loc}}

\newcommand{\D}{{\rm d}}

\newcommand{\medint}{-\kern -,375cm\int}         
\newcommand{\medintinrigo}{-\kern -,315cm\int}

\numberwithin{equation}{section}
\newtheorem{theorem}{Theorem}[section]
\newtheorem*{theorem*}{Theorem}  

\newtheorem*{expectation*}{Rivi\`ere's expectation}

\newtheorem*{conclusion*}{Conclusin}

\newtheorem*{conjecture*}{Conjecture}
\newtheorem{corollary}[theorem]{Corollary}
\newtheorem*{corollary*}{Corollary}

\newtheorem{lemma}[theorem]{Lemma}
\newtheorem*{lemma*}{Lemma}

\newtheorem*{notation*}{Notation}

\newtheorem*{problem*}{Problem}
\newtheorem{proposition}[theorem]{Proposition}
\newtheorem*{proposition*}{Proposition}
\newtheorem{remark}[theorem]{Remark}
\newtheorem*{remark*}{Remark}

\newtheorem{example}[theorem]{Example}
\newtheorem*{example*}{Example}
\newtheorem*{Acknowledgements*}{Acknowledgements}           
\theoremstyle{definition}

\makeatother

\usepackage{babel}
\begin{document}
\title[$L^p$ regularity  for even order elliptic systems]{$L^p$ regularity theory for even order elliptic systems with antisymmetric first order potentials}

 \author[C.-Y. Guo, C.-L. Xiang and G.F. Zheng ]{Chang-Yu Guo, Chang-Lin Xiang$^\ast$ and Gao-Feng Zheng}

\address[Chang-Yu Guo]{Research Center for Mathematics and Interdisciplinary Sciences, Shandong University 266237,  Qingdao, P. R. China and Institute of Mathematics, \'Ecole Polytechnique F\'ed\'erale de Lausanne (EPFL), Station 8,  CH-1015 Lausanne, Switzerland}
\email{changyu.guo@sdu.edu.cn}

\address[Chang-Lin Xiang]{Three Gorges Mathematical Research Center, China Three Gorges University,  443002, Yichang,  P. R. China, and School of Information and Mathematics, Yangtze University, 434023, Jingzhou,  P.R. China}
\email{changlin.xiang@ctgu.edu.cn}

\address[Gao-Feng Zheng]{School of Mathematics and Statistics, Central China Normal University, Wuhan 430079,  P. R.  China}
\email{gfzheng@mail.ccnu.edu.cn}


\thanks{*Corresponding author: Chang-Lin Xiang}
\thanks{C.-Y. Guo is supported by the Qilu funding of Shandong University (No. 62550089963197). The corresponding author C.-L. Xiang was financially supported by the National Natural Science Foundation of China (No. 11701045). G.-F. Zheng is supported by the National Natural Science Foundation of China (No. 11571131). 
}

\begin{abstract}	
Motivated by a challenging expectation of Rivi\`ere \cite{Riviere-2011}, in the  recent interesting  work \cite{deLongueville-Gastel-2019}, de Longueville and Gastel proposed the following geometrical even order elliptic system
\begin{equation*}
	\Delta^{m}u=\sum_{l=0}^{m-1}\Delta^{l}\left\langle V_{l},du\right\rangle +\sum_{l=0}^{m-2}\Delta^{l}\delta\left(w_{l}du\right)\qquad \text{ in } B^{2m}\label{eq: Longue-Gastel system}
\end{equation*}
which includes polyharmonic mappings as special cases. Under minimal regularity assumptions on the coefficient functions and  an additional algebraic antisymmetry assumption on the first order potential, they successfully established a conservation law for this system, from which everywhere continuity of weak solutions follows. This beautiful result amounts to a significant advance in the  expectation of Rivi\`ere.

In this paper, we seek for the optimal interior regularity of the above system, aiming at a more complete solution to the aforementioned expectation of Rivi\`ere. Combining their conservation law and some new ideas together, we obtain optimal H\"older continuity and  sharp $L^p$ regularity theory, similar to that of Sharp and Topping \cite{Sharp-Topping-2013-TAMS}, for weak solutions to a related inhomogeneous system. Our results can be applied to study heat flow and bubbling analysis for polyharmonic mappings.
\end{abstract}

\maketitle

{\small
\smallskip
\newline
\tableofcontents}
\bigskip

\section{Introduction and main results}


\subsection{Motivation}
It is well known that geometric partial differential equations are usually nonlinear in nature, for instance, the equations of harmonic mappings and the prescribed mean curvature equations. In his recent  remarkable work \cite{Riviere-2007}, to study conformally invariant variational problems in dimension two,  Rivi\`ere introduced  the following  second order linear elliptic  system \begin{equation}\label{eq:Riviere 2007}
	-\Delta u=\Omega\cdot \nabla u \qquad \text{in }B^2,
\end{equation}
where $u\in W^{1,2}(B^2, \R^n)$ and $\Omega=(\Omega_{ij})\in L^2(B^2,so_n\otimes \Lambda^1\R^2)$. As was verified in \cite{Riviere-2007}, \eqref{eq:Riviere 2007} includes
the Euler-Lagrange equations of  critical points of all second order conformally invariant variational functionals which  act on mappings $u\in W^{1,2}(B^2,N)$ from  $B^2\subset \R^2$ into a closed Riemannian manifold $N\subset \R^n$. In particular, \eqref{eq:Riviere 2007} includes the equations of weakly harmonic mappings and the prescribed mean curvature equations.

Notice that the coefficient function $\Om$ in \eqref{eq:Riviere 2007} is independent of the solution $u$. The  square integrability assumption on $\Om$ makes the system critical in the sense that $\Om\cdot \na u\in L^1(B^2)$,  which may allow discontinuous weak solutions.  Thus, from the point view of analysis, finding minimal additional assumptions on $\Om$ so that every solution would have full regularity is a very interesting and important problem. During the exploration of conformally invariant problems,  Rivi\`ere \cite{Riviere-2007} found an additional assumption which is nowadays acknowledged as the optimal one. That is, assume in addition that  $\Om $  is antisymmetric as a matrix-value function. Under this additional assumption,
Rivi\`ere found  functions $A\in L^{\wq}\cap W^{1,2}(B^2, Gl(n))$ and $B\in W^{1,2}(B^2, M_n)$
which satisfies $\nabla A-A\Omega=\na^{\bot}B$, such that system \eqref{eq:Riviere 2007} can be written  equivalently as the conservation law
\begin{equation}\label{eq:conservation law of Riviere}
	\divergence(A\nabla u+B\nabla^{\bot} u )=0,
\end{equation}
from which everywhere continuity of  weak solutions of system \eqref{eq:Riviere 2007} can be derived. As applications, this recovered the famous regularity result of  H\'elein  \cite{Helein-2002}, and confirmed affirmatively two long-standing conjectures by Hildebrandt and Heinz on  conformally invariant geometrical problems and prescribed bounded mean curvature equations respectively; see \cite{Riviere-2007} for details.

Conformally invariant variational problems  in dimensions greater than two also attracted extensive research in recent years. For fourth order conformally invariant variational problems, see for instance the  works of  Chang, Wang and Yang \cite{Chang-W-Y-1999} and  Wang \cite{Wang-2004-MZ,Wang-2004-CPAM} on biharmonic mappings from an Euclidean ball $B^n$ ($n\ge 4$) into closed Riemannian manifolds.
To extend the aforementioned powerful theory of Rivi\`ere \cite{Riviere-2007},  in the interesting work \cite{Lamm-Riviere-2008}, Lamm and Rivi\`ere proposed the following fourth order elliptic system
\begin{equation}\label{eq:Lamm-Riviere 2008}
	\De^{2}u=\De(V\cdot\na u)+{\rm div}(w\na u)+W\cdot\na u \quad  \text{in }B^4,
\end{equation}
where $V\in W^{1,2}(B^4,M_n\otimes \Lambda^1\R^{4})$, $w\in L^{2}(B^4,M_n)$,
and $W\in W^{-1,2}(B^4,M_n\otimes \Lambda^1\R^{4})$ is of the form
\[
W=\na\om+F,
\]
with $\om\in L^{2}(B^4,so_n)$ and $F\in L^{\frac{4}{3},1}(B^4,M_n\otimes \Lambda^1\R^{4})$. System \eqref{eq:Lamm-Riviere 2008} includes both extrinsic and intrinsic biharmonic mappings from $B^4$ into closed Riemannian manifolds as special cases.

Note that the first order potential $\om$ is antisymmetric as a matrix-valued function. Following the approach of  \cite{Riviere-2007},  Lamm and Rivi\`ere  \cite{Lamm-Riviere-2008} found  $A\in W^{2,2}\cap L^{\wq}(B_{1/2}^4,M_n)$ and $B\in W^{1,4/3}(B_{1/2}^4,M_n\otimes\wedge^{2}\R^{4})$
such that  $u$ is a solution of \eqref{eq:Lamm-Riviere 2008} in $B^{4}_{1/2}$ if and only if it satisfies the conservation law
\begin{equation}\label{eq:conservation law of Lamm Riviere}
	\operatorname{div}[\nabla(A \Delta u)-2 \nabla A \Delta u+\Delta A \nabla u-A w \nabla u+\nabla A(V \cdot \nabla u)-A \nabla(V \cdot \nabla u)-B \cdot \nabla u]=0.
\end{equation}
As in the second order case, everywhere continuity of weak solutions follows from  \eqref{eq:conservation law of Lamm Riviere}.

In a recent conference proceeding \cite{Riviere-2011}, Rivi\`ere gave an inspiring survey on the role of integrability by compensation in conformal geometric analysis, where in particular he presented a proof (based on his  work \cite{Riviere-2007}) on H\"older regularity of weak solutions of \eqref{eq:Riviere 2007}. Then he put up the following challenging expectation.

\begin{expectation*}\label{expect:Riviere}
	"It is natural to believe that a general result exists for $m$-th order linear systems in $m$ dimension whose 1st order potential is antisymmetric."
\end{expectation*}
\noindent A similar expectation was also  given earlier by Lamm and Rivi\`ere \cite[Remark 1.4]{Lamm-Riviere-2008}.


%




It is natural to consider conformally invariant geometrical problems in general even  dimensions, which in fact have already attracted great attention in the last decades. For instance, Gastel and Scheven \cite{Gastel-Scheven-2009CAG} obtained, among other results,  a regularity theory for both extrinsic and intrinsic polyharmonic mappings in critical dimensions. For more progress in this respect, see e.g.   \cite{Goldstein-Strzelecki-Zatorska-2009,Lamm-Wang-2009} and the references therein. Thus, a positive solution of Rivi\`ere's expectation would give, among many other applications, a new and unified approach for the H\"older regularity of (extrinsic and intrinsic) polyharmonic mappings.

However, there was no essential progress on the expectation  until quite recently. In 2019 in the interesting work \cite{deLongueville-Gastel-2019}, de Longueville and Gastel proposed   the following even order linear elliptic system
\begin{equation}\label{eq:Longue-Gastel system}
\Delta^{m}u=\sum_{l=0}^{m-1}\Delta^{l}\left\langle V_{l},du\right\rangle +\sum_{l=0}^{m-2}\Delta^{l}\delta\left(w_{l}du\right)
\end{equation}
in the unit ball $B^{2m}\subset \R^{2m}$. The coefficient functions are assumed to satisfy
\begin{equation}\label{eq:coefficient w V}
\begin{aligned}
&w_{k} \in W^{2 k+2-m, 2}\left(B^{2 m}, \mathbb{R}^{n \times n}\right) \quad \text { for } k \in\{0, \ldots, m-2\} \\
&V_{k} \in W^{2 k+1-m, 2}\left(B^{2 m}, \mathbb{R}^{n \times n} \otimes \wedge^{1} \mathbb{R}^{2 m}\right) \quad \text { for } k \in\{0, \ldots, m-1\}.
\end{aligned}
\end{equation}
Moreover,  the first order potential $V_0$ has the decomposition $V_{0}=d \eta+F$ with
\begin{equation}\label{eq:coefficient eta F}
\eta \in W^{2-m, 2}\left(B^{2 m}, s o(n)\right), \quad F \in W^{2-m, \frac{2 m}{m+1}, 1}\left(B^{2 m}, \mathbb{R}^{n \times n} \otimes \wedge^{1} \mathbb{R}^{2 m}\right).
\end{equation}
Note that $\eta$  is an antisymmetric matrix-valued function.
System \eqref{eq:Longue-Gastel system} includes both extrinsic and intrinsic $m$-polyharmonic mappings as well. It also includes equations of nonvariational type, for example ``the fake polyharmonic equations":
\[\De^m u+|\na u|^{2m} u=0, \]
provided that $u$ is bounded. For more examples, see e.g.  Strzelecki and Zatorska-Goldstein \cite{Strzelecki-Goldstein-2008-4thOrderPDE} for $m=2$.

As that of Rivi\`ere \cite{Riviere-2007} and Lamm-Rivi\`ere \cite{Lamm-Riviere-2008},   de Longueville and Gastel \cite{deLongueville-Gastel-2019}  successfully established a conservation law for system \eqref{eq:Longue-Gastel system}. More precisely,  set
\begin{equation}\label{eq:theta for small coefficient}
	\begin{aligned}
		\theta_{D}:=\sum_{k=0}^{m-2}&\|w_k\|_{W^{2k+2-m,2}(D)}+\sum_{k=1}^{m-1}\|V_k\|_{W^{2k+1-m,2}(D)}\\
		&+\|\eta\|_{W^{2-m,2}(D)}+\|F\|_{W^{2-m,\frac{2m}{m+1},1}(D)}
	\end{aligned}
\end{equation}
 for $D\subset \R^{2m}$. Under a smallness assumption
\begin{equation}\label{eq:smallness assumption}
	\theta_{B^{2m}_1}<\ep_m,
\end{equation}
they found $A\in W^{m,2}\cap L^\infty(B_{1/2}^{2m},Gl(n))$ and $B\in W^{2-m,2}(B^{2m}_{1/2},\R^{n\times n}\otimes \wedge^2\R^{2m})$ which satisfies
\begin{equation*}
\Delta^{m-1}dA+\sum_{k=0}^{m-1}(\Delta^k A)V_k-\sum_{k=0}^{m-2}(\Delta^k dA)w_k=\delta B,
\end{equation*}
such that $u$ solves \eqref{eq:Longue-Gastel system} in $B^{2m}_{{1}/{2}}$ if and only if it is a distributional solution of the conservation law
\begin{equation}\label{eq:conservation law of D-G}
\begin{aligned}
0&=\delta\Big[\sum_{l=0}^{m-1}\left(\Delta^{l} A\right) \Delta^{m-l-1} d u-\sum_{l=0}^{m-2}\left(d \Delta^{l} A\right) \Delta^{m-l-1} u \\ &\qquad -\sum_{k=0}^{m-1} \sum_{l=0}^{k-1}\left(\Delta^{l} A\right) \Delta^{k-l-1} d\left\langle V_{k}, d u\right\rangle+\sum_{k=0}^{m-1} \sum_{l=0}^{k-1}\left(d \Delta^{l} A\right) \Delta^{k-l-1}\left\langle V_{k}, d u\right\rangle \\ &\qquad -\sum_{k=0}^{m-2} \sum_{l=0}^{k-2}\left(\Delta^{l} A\right) d \Delta^{k-l-1} \delta\left(w_{k} d u\right)+\sum_{k=0}^{m-2} \sum_{l=0}^{k-2}\left(d \Delta^{l} A\right) \Delta^{k-l-1} \delta\left(w_{k} d u\right) \\ &\qquad -\langle B, d u\rangle\Big],
\end{aligned}
\end{equation}
where $d \Delta^{-1} \delta$  denotes the identity map. As an application of  \eqref{eq:conservation law of D-G}, they obtained  everywhere continuity of weak solutions of \eqref{eq:Longue-Gastel system} on $B^{2m}_{{1}/{2}}$. It is worth pointing out that in another interesting recent work of H\"orter and Lamm \cite{Horter-Lamm-2020}, the authors constructed a slightly different conservation law using a small perturbation of Uhlenbeck's gauge transform.

Before proceeding further, we would like to take a closer look at the system \eqref{eq:Longue-Gastel system} and  the  conservation law \eqref{eq:conservation law of D-G}. The first thing one may note  is the difference of coefficient functions from that of the second order system \eqref{eq:Riviere 2007} and the fourth order system \eqref{eq:Lamm-Riviere 2008}: almost half of the coefficient functions in \eqref{eq:Longue-Gastel system} are Sobolev functions with negative exponents (for definitions, see Section \ref{sec: preliminaries}). This shall cause serious problems in establishing the conservation law \eqref{eq:conservation law of D-G}. To be more accurate, to find $A,B$ as above, the authors have to solve a huge system of partial differential equations. In fact, even in the fourth order case, the task of finding $A,B$ is already quite difficult. Secondly, the  regularity of all the coefficient functions lies on the borderline of elliptic regularity theory which prevent an application of the usual $L^p$-theory, and furthermore, Sobolev functions with negative exponents require a much more general $L^p$ theory then the usual one. Due to these reasons, the authors  wrote in \cite[the last paragraph on page 19]{deLongueville-Gastel-2019} that:

\medskip
\emph{" ....... 
	But here, we consider a very general equation with rather irregular coefficients, so maybe we cannot expect much regularity in general."}
\medskip


However, there are still  several  works and results in the literature that encourage us to consider the possibility of better regularity theory than merely continuity of weak solutions to the system \eqref{eq:Longue-Gastel system}.

i) As was pointed out by  Gastel and Scheven \cite[Theorem 1.2]{Gastel-Scheven-2009CAG}, for polyharmonic mappings, H\"older continuity implies  smoothness. Thus,  H\"older continuity of solutions to system \eqref{eq:Longue-Gastel system}, if held, could be directly applied to give smoothness of weak solutions to polyharmonic mappings.

ii) The coefficients of the system \eqref{eq:Longue-Gastel system} are rather irregular, however, one observes  an additional structural feature from the system: apart from the highest order term $\De^m u$, all the left terms consist of multiplications of some derivatives of $u$, which may imply a potential iteration on the regularity of $u$. In particular, one notes that there exists no inhomogeneous term which is independent of $u$. Such an observation gives us one more support to study better regularity of weak solutions via a suitable iteration scheme.

Indeed, motivated by its geometric applications, Sharp and Topping \cite{Sharp-Topping-2013-TAMS} considered the inhomogeneous second order system (i.e., $m=1$)
 \[-\De u= \Om \cdot \na u+f\qquad \text{ in } B^2, \] and obtained   optimal H\"older continuity and sharp interior $W^{2,p}$ estimates under the assumption $f\in L^p(B^2)$ for some $1<p<2$, which  applies to Rivi\`ere's system \eqref{eq:Riviere 2007} directly. In the recent paper \cite{Guo-Xiang-Zheng-2020-Lp}, motivated also by geometric applications,  we extend the regularity result of Sharp and Topping \cite{Sharp-Topping-2013-TAMS} to the case $m=2$,  where  the  following inhomogeneous fourth order system was studied:
 \begin{equation*}
	\De^{2}u=\De(V\cdot\na u)+{\rm div}(w\na u)+W\cdot\na u+f \qquad  \text{in }B^4
\end{equation*} for some $f\in L^p(B^4)$. These two works naturally inspire us to study the regularity theory of the following inhomogeneous problem
\begin{equation}\label{eq:nonhomo Longue-Gastel system}
\Delta^{m}u=\sum_{l=0}^{m-1}\Delta^{l}\left\langle V_{l},du\right\rangle +\sum_{l=0}^{m-2}\Delta^{l}\delta\left(w_{l}du\right)+f\qquad \text{in } B^{2m}
\end{equation}
for some $f\in L^p(B^{2m})$.

iii).  $L^p$ regularity theory for \eqref{eq:nonhomo Longue-Gastel system} (when $m=1,2$) have rich geometric applications such as the energy identity (or bubbling analysis), which dates back to Sacks-Uhlenbeck \cite{Sacks-Uhlenbeck-1981} on harmonic mappings, and heat flow which dates back to e.g. Struwe \cite{Struwe-1985,Struwe-1988}. Indeed, the work of Sharp and Topping \cite{Sharp-Topping-2013-TAMS} has many interesting applications in these problems. For instance, it has been successfully applied in the study of angular energy quantization of weak solutions of system \eqref{eq:Riviere 2007} by Laurain and Rivi\`ere \cite{Laurain-Riviere-2014-APDE}, which largely extended the bubbling analysis of Sacks-Uhlenbeck \cite{Sacks-Uhlenbeck-1981}. In the fourth order case, $L^p$ estimates have also been applied to bubbling analysis  in  e.g. \cite{Laurain-Riviere-2013-ACV,Hornung-Moser-2012,Wang-Zheng-2012-JFA}.  $L^p$ regularity estimates have also found applications in global estimates, which leads to the full energy quantization results; see for instance the interesting work and Lamm and Sharp \cite{Lamm-Sharp-2016-CPDE}. For extension and applications to supercritical dimensional problems, see e.g. \cite{Sharp-2014,Moser-2015-TAMS}.

Based on the above considerations and seeking for a more complete regularity theory for Rivi\`ere's expectation, it is natural to  ask
\begin{problem*}\label{quest:main question}
	What is the optimal H\"older continuity and $L^p$ regularity for  weak solutions of \eqref{eq:nonhomo Longue-Gastel system}? Can we derive optimal interior estimates for  \eqref{eq:nonhomo Longue-Gastel system} as that of \cite{Sharp-Topping-2013-TAMS} for $m=1$ and \cite{Guo-Xiang-Zheng-2020-Lp} for $m=2$?
\end{problem*}

The aim of this paper is to give an affirmative answer to this problem, and thus, provide a complete solution to Rivi\`ere's  expectation in the even order case.

\subsection{Main results}

Our first theorem deals with the optimal H\"older continuity of weak solutions to  \eqref{eq:Longue-Gastel system}.

\begin{theorem}[H\"older continuity]\label{thm:optimal Holder exponent for inho DG}
Suppose $f\in L^p(B^{2m}_{1})$ for some $p\in (1,\frac{2m}{2m-1})$. If $u\in W^{m,2}(B^{2m}_{1},\R^n)$ is a weak solution of  \eqref{eq:nonhomo Longue-Gastel system}, then $$u\in C^{0,\al}_{\loc}(B^{2m}_{1})$$ with  $\al=2m(1-{1}/{p})$.
Moreover, there exist  $r_0, C>0$ depending only on $m,n,p$ and the coefficient functions $V_k,w_k$ such that for all $0<r<r_0$, there holds
 \begin{equation}\label{eq: decay estimate for 2m order}
    \begin{aligned}
    \sum_{j=1}^m\|\na^j u\|_{L^{2m/j,2}(B^{2m}_{r})}\leq  Cr^{\al}\left( \sum_{j=1}^m\|\na u\|_{L^{2m/j,2}(B^{2m}_1)}+\|f\|_{L^p(B^{2m}_{1})}\right).
      \end{aligned}
     \end{equation}
\end{theorem}

The H\"older continuity is optimal, as one can see from the simplest case $\De^m u=f$ by noticing the Sobolev embedding $W^{2m,p}(B^{2m}_{1})\subset C^{0,2m(1-1/p)} (B^{2m}_{1})$. Moreover,  applying Theorem \ref{thm:optimal Holder exponent for inho DG}  to the case $f\equiv 0$ yields
 that every weak solution $u\in W^{m,2}(B^{2m}_{1})$ of system \eqref{eq:Longue-Gastel system} is locally $\al$-H\"older continuous for all $\al\in (0,1)$\footnote{In \cite{Guo-Xiang-2019-Higher} we proved  $\al$-H\"older continuity  for solutions to \eqref{eq:Longue-Gastel system}  for some $\al\in (0,1)$ without using  conservation law.}. On the other hand, the H\"older continuity is the best possible regularity that one can expect for weak solutions of \eqref{eq:Longue-Gastel system}. In the case $m=2$, we constructed an example in \cite{Guo-Xiang-Zheng-2020-Lp} which fails to be locally Lipschitz continuous. Using a similar idea, it is not difficulty to construct such an example for the general case and we leave it to the interested readers.

In our second theorem, we derive optimal higher order regularity of weak solutions.

\begin{theorem}[local $L^p$ estimates]\label{thm:optimal global estimate for inho DG}
Let $u\in W^{m,2}(B^{2m}_{1},\R^n)$ be a weak solution of \eqref{eq:nonhomo Longue-Gastel system} with $f\in L^p(B^{2m}_{1})$ for $p\in (1,\frac{2m}{2m-1})$. Then $$u\in W_{\loc}^{m+1,\frac{2mp}{2m-(m-1)p}}(B^{2m}_{1}).$$ Moreover,
 there exist $\ep>0$ and $C>0$ depending only on $m,n,p$ and the coefficient functions $V_k,w_k$  such that,  if the smallness condition \eqref{eq:smallness assumption} is satisfied on $B^{2m}_{1}$ with $\ep_m=\ep$, then
    \begin{equation}\label{eq:optimal m+1 order}
    \|u\|_{W^{m+1,\frac{2mp}{2m-(m-1)p}}(B^{2m}_{1/2})}\le C\left(\|f\|_{L^{p}(B^{2m}_{1})}+\|u\|_{L^{1}(B^{2m}_{1})}\right).
    \end{equation}
\end{theorem}

Theorem \ref{thm:optimal global estimate for inho DG} can be regarded as a counterpart of Sharp-Topping \cite[Theorem 1.1]{Sharp-Topping-2013-TAMS} (corresponding to the case $m=1$) and our work \cite[Theorem 1.2]{Guo-Xiang-Zheng-2020-Lp} (corresponding to the case $m=2$) to the general even order system \eqref{eq:Longue-Gastel system}. The order $m+1$ is also the best possible; see Example \ref{example:no Wm+2 estimate} below (and the case $m=2$ was explained in \cite{Guo-Xiang-Zheng-2020-Lp}).

 As a direct  application of Theorem \ref{thm:optimal global estimate for inho DG}, we have the following improvement of the regularity result of de Longueville and Gastel \cite{deLongueville-Gastel-2019}.
\begin{corollary}
  If $u\in W^{m,2}(B_{1},\R^n)$ is a weak solution of  \eqref{eq:Longue-Gastel system}, then
  \[u\in \bigcap_{0<\ep<1}W^{m+1,2-\ep}_{\loc}(B^{2m}_1)\subset  \bigcap_{0<\al<1} C^{0, \al}_{\loc}(B^{2m}_1).\]
\end{corollary}

One more application of Theorem \ref{thm:optimal global estimate for inho DG} yields the following energy gap, which was known to be useful in deducing the energy quantization results; see for instance  \cite{Lin-Riviere-2002,Laurain-Riviere-2014-APDE,Hornung-Moser-2012}.

\begin{corollary}[Energy gap]\label{coro:energy gap}
	Let $u\in W^{m,2}(\R^{2m},\R^n)$ be a weak solution of \eqref{eq:Longue-Gastel system} in $\R^{2m}$.  There exists some $\ep=\ep(m,n)>0$ such that if
$$\theta_{\R^{2m}}<\ep,$$
 then $u\equiv 0$ in $\R^{2m}$, where $\theta_{\R^{2m}}$ is defined as in \eqref{eq:theta for small coefficient}.
\end{corollary}

To control the size of the paper, we do not consider in this paper the compactness problem of weak solutions of system \eqref{eq:nonhomo Longue-Gastel system} under the weak $L\log L$-integrability of $f$, as that of Sharp-Topping \cite{Sharp-Topping-2013-TAMS} (when $m=1$) and Guo-Xiang-Zheng \cite{Guo-Xiang-Zheng-2020-Lp} (when $m=2$).

\subsection{Strategy of the proof}

In spirit, we follow the scheme of Sharp and Topping \cite{Sharp-Topping-2013-TAMS} and \cite{Guo-Xiang-Zheng-2020-Lp}  to divide  the proof into three main steps.
\begin{enumerate}
	\item In the first step, we  derive the Morrey type decay estimate in Theorem \ref{thm:optimal Holder exponent for inho DG} via the conservation law \eqref{eq:conservation law of D-G}, from which H\"older continuity follows.  The decay estimates show that $\na^j u$ ($1\leq j\leq m$) belong to some Morrey spaces.
	
	\item In the second step, we combine the above Morrey type estimate, together with the Riesz potential theory of Adams \cite{Adams-1975} (see Lemma \ref{lemma:improved Riesz potential}) on Morrey spaces, to deduce an almost optimal higher order Sobolev regularity. That is, we prove that $u\in W^{m+1, q}_{\loc}$ for all $q$ smaller than the optimal exponent $\bar{q}=\frac {2mp}{2m-(m-1)p}$.
	
	\item In the final step, we show that  the concerned  local $W^{m+1, q}$ estimates are uniform with respect to $q<\bar{q}$ so that by passing to the limit we can obtain the asserted optimal higher order Sobolev regularity.
\end{enumerate}

As we are dealing with a very general even order system with highly irregular coefficients, it is expected that the situation will become much more complicated than that of the second order system considered in Sharp and Topping \cite{Sharp-Topping-2013-TAMS},  where several parts require more delicate treatments.

The first step is fundamental for the following steps and it is based on an iteration procedure. In Sharp and Topping \cite{Sharp-Topping-2013-TAMS}, their smart iteration argument  requires a very precise control on some coefficients of related estimates so that the key coefficients are sufficiently small. Such a precise control is based on the monotonicity of the function $r\mapsto |B_r|^{-1}\int_{B_r }f$ for subharmonic functions $f$. Thus it fails even in the fourth order case $m=2$.    In our case, it is far more subtle and a more general treatment is required: we have to use an extremely complicated conservation law \eqref{eq:conservation law of D-G} to figure out the system for $A\Delta u$ and $A\nabla u$. As there are many terms with low regularity appearing, we have to do a very careful case study in order to properly apply the Riesz potential theory to gain the desired decay estimates.

Another significant difference occurs in the last step. In the case of Sharp-Topping, they run a similar scaling and iteration argument as in the first step due to the well control of coefficient of decay of energy of $\nabla u$. Moreover,  in the second order case every solution belongs to $W^{2,p}_{\loc}$, which allows them to view the system $-\De u=\Om\cdot \na u$  as a pointwise identity from that they can deduce  estimate for $\na^2 u$ directly. In our case, it is far more difficult to obtain a uniform estimate for  $\|\nabla^{m+1}u\|_q$. To overcome this difficulty, we  combine a duality argument, together with some elliptic estimates on polyharmonic maps with Navier boundary condition, to obtain an estimate such that Simon's iteration method can be applied to deduce a quantitative and uniform control on the norm $\|\nabla^{m+1}u\|_{q, B_{1/2}}$ for all $q$ smaller than the optimal exponent. A disadvantage of this argument is that  we only  obtain an upper bound of the form $\|u\|_{W^{m,2}}+\|f\|_{L^p}$. To improve this bound,  we will have to run an additional interpolation argument to obtain  the better bound in \eqref{eq:optimal m+1 order}.

This paper is organized as follows. Section \ref{sec: preliminaries} contains some preliminaries and auxiliary results for later proofs. Our main theorems are proved in Sections \ref{sec:Decay estimate forth system}, \ref{sec:higher order Sobolev regularity} and \ref{sec:optimal global estimates}. We also add an appendix to include certain auxiliary results that was used in the proofs of our main theorems.

Our notations are standard.
By $A\lesssim B$ we mean there exists a universal constant $C>0$ such that $A\le CB$.

\begin{Acknowledgements*}
  We are grateful to professor Xiao Zhong for  many useful suggestions and discussions during the preparation of this manuscript.
\end{Acknowledgements*}

\section{Preliminaries and auxiliary results}\label{sec: preliminaries}

\subsection{Function spaces and related}\label{sec:relevant function spaces}
Let $\Om\subset\R^{n}$ be a bounded smooth  domain, $1\le p<\wq$
and $0\le s<n$. The \emph{Morrey space} $M^{p,s}(\Om)$ consists
of functions $f\in L^{p}(\Om)$ such that
\[
\|f\|_{M^{p,s}(\Om)}\equiv\left(\sup_{x\in\Om,r>0}r^{-s}\int_{B_{r}(x)\cap\Om}|f|^{p}\right)^{1/p}<\wq.
\]
Denote by $L_*^p$  the weak $L^p$ space and define the \emph{weak Morrey space} $M^{p,s}_*(\Omega)$ as the space of functions $f\in L^p_*(\Omega)$ such that
$$\|f\|_{M^{p,s}_*(\Om)}\equiv \left(\sup_{x\in\Om,r>0}r^{-s}\|f\|^p_{L^p_{\ast}(B_r(x)\cap \Omega)}\right)^{1/p}<\wq,$$
where
$$\|f\|^p_{L^p_{\ast}(B_r(x)\cap \Omega)}\equiv\sup_{t>0}t^p\Big|\left\{x\in B_r(x)\cap \Omega: |f(x)|> t\right\}\Big|.$$

For a measurable function $f\colon \Om\to\R$, denote by $\de_{f}(t)=|\{x\in\Om:|f(x)|>t\}|$
its distributional function and by $f^{\ast}(t)=\inf\{s>0:\de_{f}(s)\le t\}$,
$t\ge0$, the nonincreasing rearrangement of $|f|$. Define
\begin{eqnarray*}
	f^{\ast\ast}(t)\equiv\frac{1}{t}\int_{0}^{t}f^{\ast}(s)\D s, &  & t>0.
\end{eqnarray*}
The\emph{ Lorentz space} $L^{p,q}(\Om)$ ($1<p<\wq,1\le q\le\wq$)
is the space of measurable functions $f:\Om\to\R$ such that
\[
\|f\|_{L^{p,q}(\Om)}\equiv\begin{cases}
\left(\int_{0}^{\wq}(t^{1/p}f^{\ast\ast}(t))^{q}\frac{\D t}{t}\right)^{1/q}, & \text{if }1\le q<\wq,\\
\sup_{t>0}t^{1/p}f^{\ast\ast}(t) & \text{if }q=\wq
\end{cases}
\]
is finite.

  It is well-known that $L^{p,p}=L^p$ and $L^{p,\wq}=L^p_{\ast}.$ We will need the following H\"older's inequality in Lorentz spaces.

\begin{proposition}[\cite{ONeil-1963}]\label{prop: Lorentz-Holder inequality}
	Let $1<p_{1},p_{2}<\wq$ and $1\le q_{1},q_{2}\le\wq$ be such that
	\begin{eqnarray*}
		\frac{1}{p}=\frac{1}{p_{1}}+\frac{1}{p_{2}}\le1 & \text{and} & \frac{1}{q}=\frac{1}{q_{1}}+\frac{1}{q_{2}}\le1.
	\end{eqnarray*}
	Then, $f\in L^{p_{1},q_{1}}(\Om)$ and $g\in L^{p_{2},q_{2}}(\Om)$
	implies $fg\in L^{p,q}(\Om)$. Moreover,
	\[
	\|fg\|_{L^{p,q}(\Om)}\le\|f\|_{L^{p_{1},q_{1}}(\Om)}\|g\|_{L^{p_{2},q_{2}}(\Om)}.
	\]
\end{proposition}

\begin{proposition}[\cite{Ziemer}]\label{prop: Lorentz to Lorentz}
	For $1<p<\wq$ and $1\le q_{1}\le q_{2}\le\wq$, there holds $L^{p}(\Om)=L^{p,p}(\Om)$
	and $L^{p,q_{1}}(\Om)\subset L^{p,q_{2}}(\Om)$ with
	\[
	\|f\|_{L^{p,q_{2}}(\Om)}\le C(p,q_1,q_2)\|f\|_{L^{p,q_{1}}(\Om)}.
	\]
	Moreover, if $|\Om|<\wq$, then $L^{p,q}(\Om)\supset L^{r,s}(\Om)$
	for $1<p<r<\wq$, $1\le q,s\le\wq$, and
	\[
	\|f\|_{L^{p,p}(\Om)}\le C_{r,p} |\Om|^{\frac{1}{p}-\frac{1}{r}}\|f\|_{L^{r,\wq}(\Om)}.
	\]
\end{proposition}




Lorentz spaces can be used to define the following Lorentz-Sobolev spaces.
For $k\in \mathbb{N}$, $1<p<\wq, 1\leq q\leq \infty$, the Lorentz-Sobolev space $W^{k,p,q}(\Omega)$ consists  of  measurable functions $f\colon \Omega\to \R$ which are weakly differentiable up to order $k$ with $\na^{\alpha}f\in L^{p,q}(\Omega)$ for all $|\alpha|\leq k$, and its norm is defined by
\[
\|f\|_{W^{k,p,q}(\Om)}=\left(\sum_{|\al|\le k}\|\na^{\al} f\|_{L^{p,q}(\Om)}^{p}\right)^{1/p}.\]

\begin{proposition}[\cite{Tartar-1998}]\label{prop:Lorentz-Sobolev embedding}
	Let $k\in\N, 1<p<\wq$ and $1\le q\le \wq$.  Then
	
	\begin{itemize}
		\item[(1).] $W^{k,p}(\Om)=W^{k,p,p}(\Om)$;
		\item[(2).] If $\Omega$ is bounded and smooth and $kp<n$, then $W^{k,p,q}(\Om)$ embeds continuously into $L^{p^{*},q}(\Om)$,
	where $1/p^{*}=1/p-k/n$.
         \item[(3).]  $W^{k,n/k,1}(\Om)\subset C(\overline{\Om})$.
	\end{itemize}
\end{proposition}

\begin{proposition}\label{prop: extension thm} (\cite[Theorem 2]{DeVore-Scherer-1979}, \cite[Page 19]{deLongueville-2018})
   Let $\Om\subset\R^n$ be a bounded smooth domain and $V\subset\R^n$ an open set such that $\Om\subset\subset V$, $k\in \N , 1<p<\wq, 1\le q\le \wq$. There exists a bounded linear operator $E: W^{k,p,q}(\Omega)\to W^{k,p,q}(\R^n)$, such that for every  $f\in W^{k,p,q}(\Omega)$,

    (i) $Ef=f$ a.e. in $\Om$,

     (ii) $Ef$ is compactly supported in $V$. Moreover, there exists a constant $C>0$ depending only on $n,k,p$ and $\Om$,  such that for all $f\in W^{k,p,q}(\Omega)$, there holds
   \[\|Ef\|_{W^{k,p,q}(\R^n)}\le C\|f\|_{W^{k,p,q}(\Omega)}.\]
\end{proposition}

For later use, we furthermore need  Lorentz-Sobolev spaces with negative exponents. For $k\in \N$, $1\leq p,q\leq \infty$, the Lorentz-Sobolev space $W^{-k,p,q}(\Omega)$ is defined to be the space of all distributions $f$ on $\Om$ of the form $f=\sum_{|\alpha|\leq k}\na^{\alpha}f_{\alpha}$ with $f_\alpha\in L^{p,q}(\Om)$. The corresponding norm is defined as
\[ \|f\|_{W^{k,p,q}(\Omega)}:=\inf \sum_{|\alpha|\leq k}\|f_\alpha\|_{L^{p,q}(\Omega)},\]
where the infimum is taken over all decompositions of $f$ as given in the definition.
 The following facts on Lorentz-Sobolev spaces with negative exponents can be found in \cite{deLongueville-2018,deLongueville-Gastel-2019}.
\begin{proposition}\label{prop: generalized holder ineq}
Let $\Omega\subset \R^n$ be a bounded smooth domain.

\begin{itemize}
\item[(1)] {\rm (Generalized H\"older's inequality)} Suppose $f\in W^{-k_0,p_0,q_0}(\Omega)$ and $g\in W^{k_1,p_1,q_1}(\Omega)$ with $k_0,k_1\in \N$, $1<p_0,p_1<\infty$ and $1\leq q_0,q_1<\infty$. If $k_0\le k_1$, $\frac{1}{p_0}+\frac{1}{p_1}\leq 1$ and $k_1p_1<n$, then $fg\in W^{-k_0,s,t}(\Omega)$ with $s=\frac{np_0p_1}{n(p_0+p_1)-k_1p_0p_1}$ and $\frac{1}{t}=\min\big\{1,\frac{1}{q_0}+\frac{1}{q_1}\big\}$. Moreover,
\begin{equation}\label{eq:product norm for Lorentz-Sobolev}
\|fg\|_{W^{-k_0,s,t}(\Om)}\leq C\|f\|_{W^{-k_0,p_0,q_0}(\Om)}\|g\|_{W^{k_1,p_1,q_1}(\Om)}.
\end{equation}
The assertion continues to hold in the case $k_1p_1=n$ if we additionally assume $g\in L^\wq$.

\item[(2)]{\rm (Sobolev embedding)} If $k\in \mathbb{Z}$, $l\in \N$, $1<p<\frac{n}{l}$ and $1\leq q\leq \infty$, then $W^{k,p,q}(\Omega)$ embeds continuously into $W^{k-l,\frac{np}{n-lp},q}(\Omega)$ with
$$\|f\|_{W^{k-l,\frac{np}{n-lp},q}(\Omega)}\leq C\|f\|_{W^{k,p,q}(\Omega)}.$$
\end{itemize}
\end{proposition}

\subsection{Fractional Riesz operators}

Let $0<\al<n$ and  $I_\alpha=c_{n,\al}|x|^{\alpha-n}$, $x\in \R^n$, be the usual fractional Riesz operators, where $c_{n,\al}$ is a positive normalization constant. The following well-known estimates on fractional Riesz operators in Morrey spaces were proved by Adams \cite{Adams-1975}.

\begin{proposition}\label{prop:Adams 1975} Let $0< \alpha<n$, $0\le \lambda< n$ and $1\le p<\frac{n-\lambda}{\alpha}$. There exists a constant $C>0$ depending only $n,\al,\la$ and $p$ such that, for all $f\in M^{p,\lambda}(\R^n)$, there holds

	(i) If  $p>1$, then
 \begin{equation}\label{eq:Riesz Adams 1}
		\|I_\alpha(f)\|_{M^{\frac{(n-\la)p}{n-\la-\al p},\lambda}(\R^n)}\leq C\|f\|_{M^{p,\lambda}(\R^n)}.
		\end{equation}
					
		(ii) If $p=1$, then
		\begin{equation}\label{eq:Riesz Adams 2}
		\|I_\alpha(f)\|_{M_{*}^{\frac{n-\lambda}{n-\lambda-\alpha},\lambda}(\R^n)}\leq C\|f\|_{M^{1,\lambda}(\R^n)}.
		\end{equation}	
\end{proposition}
Note that when $\la=0$ it reduces to  the usual   Riesz potential theories between $L^p$ spaces.

\begin{proposition}\label{prop:Riesz potential}
For $0<\alpha<n$, $1<p<n/\al$, $1\le q\leq q'\le\wq$, the fractional Riesz operators
\[
I_{\al}\colon L^{p,q}(\R^{n})\to L^{\frac{np}{n-\al p},q'}(\R^{n})
\]
and
\[
I_{\al}\colon L^{1}(\R^{n})\to L^{\frac{n}{n-\al},\wq}(\R^{n})
\]
are bounded.
\end{proposition}
In the later proof, we will mainly use the special cases
\[
I_{\al}:L^{\frac{2m}{\al+j},2}(\R^{2m})\to L^{\frac{2m}{j},2}(\R^{2m})
\]
for $0<\al<2m$ and $m\ge j\ge1$ such that $\al+j<2m$.

The  result of Proposition \ref{prop:Riesz potential} continues to hold for $\al=0$ (see Theorems V.3.15 and VI.3.1 of Stein and Weiss \cite{Stein-Weiss}), provided we admit to use the notation $I_0$ to denote singular integral operators of the type $\na^{2m-k}I_k$ for  integers $1\le k\le 2m-1$. In these cases, it will be convenient to introduce the notation $I_0$.




The following lemma plays a crucial role in the proofs of Sharp-Topping \cite{Sharp-Topping-2013-TAMS}, and \cite{Guo-Xiang-Zheng-2020-Lp}, and  also in our approach here. It can be viewed as an improved version of the classical Riesz potential theory.

\begin{lemma}[\cite{Adams-1975}, Proposition 3.1]\label{lemma:improved Riesz potential}
	Let $0<\al<\be\le n$ and $f\in M^{1,n-\be}(\R^{n})\cap L^{p}(\R^{n})$
	for some $1<p<\wq$. Then, $I_{\al}f\in L^{\frac{\be p}{\be-\al}}(\R^{n})$
	with
	\[
	\|I_{\al}f\|_{\frac{p\be}{\be-\al}, \R^n}\le C_{\al,\be,n,p}\|f\|_{M^{1,n-\be}(\R^n)}^{\frac{\al}{\be}}\|f\|_{p, \R^n}^{\frac{\be-\al}{\be}}. 
	\]
	\end{lemma}
In particular, in the case $1<p<n/\be$, we have
	\[
	\frac{\be p}{\be-\al}>\frac{np}{n-\al p},
	\]
	which implies that $I_{\al}f$ has better integrability than the usual one from $L^p$ boundedness. This is due to the fact that $f$ has additional fine property, that is, $f$ also belongs to some
Morrey space. For a local version of Lemma \ref{lemma:improved Riesz potential}, see  \cite[Lemma A.3]{Sharp-Topping-2013-TAMS}.

%
%
%

\subsection{Scaling invariance of \eqref{eq:nonhomo Longue-Gastel system}}\label{sec: scaling invariance}
We shall  use a scaling argument in our later proofs. Let $u$ be a weak solution of  \eqref{eq:nonhomo Longue-Gastel system}.  For any  $B_{2R}(x_{0})\subset B_{1}\subset \R^{2m}$ and $x\in B_{1/2}=B_{1/2}(0)$, set
\[
\begin{aligned} & u_{R}(x)=u(x_{0}+Rx),\quad V_{l,R}(x)=R^{2m-2l-1}V_l(x_{0}+Rx),\\
	&w_{l,R}(x)=R^{2m-2l-2}w_l(x_{0}+Rx),\quad f_{R}(x)=R^{2m}f(x_{0}+Rx).
\end{aligned}
\]
It is straightforward to verify that $u_R$ satisfies
\begin{eqnarray}\label{eq: scaled equation}
\Delta^{m}u_R=\sum_{l=0}^{m-1}\Delta^{l}\left\langle V_{l,R},du_R\right\rangle +\sum_{l=0}^{m-2}\Delta^{l}\delta\left(w_{l,R}du_R\right)+f_R &  & \text{in }B_{1/2}.
\end{eqnarray}
Moreover, for any $0<r<1$ and $1\le q\le \wq$,  there hold
\[
\begin{aligned}
	\|\na^i u_{R}\|_{L^{\frac{2m}{i},q}(B_{r}(0))}=\|\na^i u\|_{L^{\frac{2m}{i},q}(B_{rR}(x_{0}))},\quad\|f_{R}\|_{L^{p}(B_{r}(0))}=R^{2m(1-1/p)}\|f\|_{L^{p}(B_{rR}(x_{0}))}.
\end{aligned}
\]

\section{H\"older continuity via decay estimates}\label{sec:Decay estimate forth system}
In this section, we prove  Theorem \ref{thm:optimal Holder exponent for inho DG}.
Throughout this section, we use $B_{r}(x)$ to denote the ball with centre $x$ and radius $r$ in $\R^{2m}$, and    write $B_r$ instead of $B_r(0)$ when $x=0$. Set
 \[
1<p<2m/(2m-1)\qquad \text{and}\qquad \al_{0}=2m(1-1/p).
\]
To prove Theorem \ref{thm:optimal Holder exponent for inho DG}, it suffices to establish the following decay estimate.

 \begin{lemma}[Decay estimate]\label{lem: decay of Long-Gastel}
Suppose $f\in L^p(B_1)$ for some $1<p<2m/(2m-1)$ and $u\in W^{m,2}(B_{1}, \R^{n})$ is a solution to system \eqref{eq:nonhomo Longue-Gastel system}. There exist $r_{0},\tau\in(0,1)$ and $C>0$ depending
only on $m,n,p$ and the coefficient functions $V_k, w_k$ such that for $x\in B_{1/2}$
and $0<r<r_{0}$,
\[
\sum_{j=1}^{m}\|\na^{j}u\|_{L^{2m/j,2}(B_{r}(x))}\le C\left(\sum_{j=1}^{m}\|\na^{j}u\|_{L^{2m/j,2}(B_{1})}+\|f\|_{L^{p}(B_{1})}\right)r^{\al_{0}}.
\]
\end{lemma}

Once the above lemma is proved, Theorem \ref{thm:optimal Holder exponent for inho DG} follows easily.

\begin{proof}[Proof of Theorem \ref{thm:optimal Holder exponent for inho DG}] By Lemma \ref{lem: decay of Long-Gastel} and Proposition \ref{prop: Lorentz to Lorentz}, we have
\[
\int_{B_{r}(x)}|\na u|^{2m}\le C\|\na u\|_{L^{2m,2}(B_{r}(x))}^{2m}\le C\left(\|u\|_{W^{m,2}(B_{1})}+\|f\|_{L^p(B_1)}\right)^{2m}r^{2m\al_0}
\]
for all $x\in B_{1/2}$ and  $0<r<r_{0}$. Then Morrey's Dirichlet
growth theorem \cite{Giaquinta-Book} implies that $u\in C^{0,\al_0}(B_{1/2})$. The proof
is complete. \end{proof}

Due to the scaling invariance exhibited in section \ref{sec: scaling invariance},  Lemma \ref{lem: decay of Long-Gastel} is a consequence of the following estimate.
\begin{lemma}\label{lem: decay estimate lemma 2}
 Suppose $f\in L^p(B_1)$ for some $1<p<2m/(2m-1)$ and $u\in W^{m,2}(B_{1}, \R^{n})$ is a solution to system \eqref{eq:nonhomo Longue-Gastel system}.  There exist $\ep,C>0$ depending only on $m,n,p$ and the coefficient functions $V_k, w_k$ such that, if
\[  \ta_{B_1}<\ep,\] then
\begin{equation}
\sum_{j=1}^{m}\|\na^{j}u\|_{L^{2k/j,2}(B_{\tau})}
\le\tau^{\be}\sum_{j=1}^{m}\|\na^{j}u\|_{L^{2k/j,2}(B_1)}
+C\tau^{\al_{0}}\|f\|_{L^{p}(B_{1})},
\label{eq: decay estimates of higher order system}
\end{equation} where $\ta_{B_1}$ is defined in \eqref{eq:theta for small coefficient} and $\be=(\al_{0}+1)/2$.
\end{lemma}

The proof of Lemma \ref{lem: decay estimate lemma 2} is lengthy.  We divide it into four steps.   First  choose a sufficiently small $\epsilon$ in Lemma \ref{lem: decay estimate lemma 2} such that the conservation law \eqref{eq:conservation law of D-G} holds for some $A,B$ in $B_1$\footnote{Recall that  under the smallness assumption $\ta_{B_1}<\ep$ the conservation law \eqref{eq:conservation law of D-G}  only holds on $B_{1/2}$, but here for simplicity  we directly assume it hold on the whole ball $B_1$ since otherwise we can do everything below on  $B_{1/2}$}.

\textbf{Step 1.} Deduce the equation of $Adu$.
\begin{proposition}\label{prop:equation for A nabla u}
$Adu$ satisfies the system
\begin{equation}
\de\De^{m-1}(Adu)=\sum_{i=1}^{m-1}\de^{i}\left(\sum_{j=m-i}^{m}\na^{j}A\na^{2m-i-j}u\right)+\de K+Af,\label{eq: equation of Adu}
\end{equation}
where $K$ is the last five terms of the conservation law \eqref{eq:conservation law of D-G}.
\end{proposition}

In below $\de$ denotes the divergence operator and $\de^i$  means taking divergence for $i$ times. By $\na^j A\na^{k}u$ we mean it is a linear combination of products of the type  $\na^{\al} A \na^{\be} u$ for which $|\al|=j$ and $|\be|=k$.

\begin{proof} Repeatedly using Leibniz rule gives
\begin{equation}
\de\De^{m-1}(Adu)=\de(A\De^{m-1}du)
+\sum_{i=1}^{m-1}\de^{i}\left(\sum_{j=m-i}^{m}\na^{j}A\na^{2m-i-j}u\right).\label{eq: 1}\footnote{A proof of this equality can be found in the Appendix.}
\end{equation}
Then the conservation law \eqref{eq:conservation law of D-G} yields
\[
\de(A\De^{m-1}du)=\de\left(\sum_{l=1}^{m-1}\De^{l}A\De^{m-1-l}du+\sum_{l=0}^{m-2}d\De^{l}A\De^{m-1-l}u\right)+\de K+Af,
\]
where $K$ is the last five terms of the conservation law of \eqref{eq:conservation law of D-G}. Hereafter, we use $\sum_i a_i$ to denote a linear combination of $a_i$'s, i.e., $\sum_i a_i= \sum_i c_i a_i$  for some harmless absolute constants $c_i$. We will  write down the coefficient explicitly when necessary.

Furthermore, note that the first term in the right hand side of the above equality can be written as
\[
\de\left(\sum_{l=1}^{m-1}\De^{l}A\De^{m-1-l}du+
\sum_{l=0}^{m-2}d\De^{l}A\De^{m-1-l}u\right)
=\sum_{i=1}^{m-1}\de^{i}\left(\sum_{j=m-i}^{m}\na^{j}A\na^{2m-i-j}u\right).
\]
 For instance, repeatedly using Leibniz rule gives (with $l=1$)
\[
\De A\De^{m-2}du=\sum_{i=0}^{m-2}\de^{i}(\na^{m-i}A\na^{m-1}u),
\]
and (with $l=m-1$)
\[
\De^{m-1}Adu=\sum_{i=0}^{m-2}\de^{i}(\na^{m}A\na^{m-i}u),
\]
and (with $l=0$)
\[
dA\De^{m-1}u=dA\De^{m-2}\De u=\sum_{i=0}^{m-2}\de^{i}(\na^{m-1-i}A\na^{m}u)
\]
and (with $l=m-2$)
\[
d\De^{m-2}A\De u=\De^{m-2}dA\De u=\sum_{i=0}^{m-2}\de^{i}(\na^{m-1}A\na^{m-i}u).
\]
This leads to \eqref{eq: equation of Adu}. The proof is complete.
\end{proof}

By noting that
\[
\De^{m-1}(A\De u)=\de\De^{m-1}(Adu)-\De^{m-1}(dAdu),
\]
we infer  from Proposition \ref{prop:equation for A nabla u} that
\begin{corollary}\label{coro:equation for A Delta u}
$A\Delta u$ satisfies system
\begin{equation}\label{eq:A Delta u}
\De^{m-1}(A\De u)=\sum_{i=1}^{m-1}\de^{i}\left(\sum_{j=m-i}^{m}\na^{j}A\na^{2m-i-j}u\right)+\de K+Af,
\end{equation}
where $K$ is the last five terms of the conservation law \eqref{eq:conservation law of D-G}.
\end{corollary}

\textbf{Step 2.} Decomposition of $Adu$.

We use  Hodge decomposition to obtain
\begin{eqnarray*}
Adu=d\tilde{u}_{1}+d^{\ast}\tilde{u}_{2}+\tilde{h} &  & \text{in }B_{1},
\end{eqnarray*}
such that
\[
\De^{m}\tilde{u}_{1}=\sum_{i=1}^{m-1}\de^{i}\left(\sum_{j=m-i}^{m}\na^{j}A\na^{2m-i-j}u\right)+\de K+Af,
\]
\[
\De^{m}\tilde{u}_{2}=\De^{m-1}(dA\wedge du)
\]
in $B_{1}$ and $\tilde{h}$ is a harmonic one form in $B_{1}$. Here to be more precise, we use $d^{\ast}$ to mean the codifferential operator.

Now we extend all the functions from $B_{1}$ into the whole space
$\R^{2m}$ in a bounded way. Still we use the same notations to denote
the extended functions.
We define
\begin{equation}
u_{1}=c\log\ast\left(\sum_{i=1}^{m-1}\de^{i}\left(\sum_{j=m-i}^{m}\na^{j}A\na^{2m-i-j}u\right)+\de K+Af\right)\label{eq: singular part}
\end{equation}
and
\begin{equation}
u_{2}=c\log\ast\De^{m-1}(dA\wedge du)\label{eq: def of u2}
\end{equation}
such that $\De^{m}(\tilde{u}_{1}-u_{1})=\De^{m}(\tilde{u}_{2}-u_{2})=0$
in $B_{1}$, where $c\log$ is the fundamental solution of $\De^m$ in $\R^{2m}$, and  then define $h=\tilde{u}_{1}-u_{1}+\tilde{u}_{2}-u_{2}+\tilde{h}$.
Consequently, we obtain
\begin{eqnarray}
Adu=du_{1}+ d^{\ast}u_{2}+h &  & \text{in }B_{1},\label{eq: hodge decomposition}
\end{eqnarray}
 with $u_{1}$ given by (\ref{eq: singular part}), $u_{2}$ given
by (\ref{eq: def of u2}), and $h$ is an $m$-polyharmonic 1-form
in $B_{1}$, i.e.,
\begin{eqnarray*}
\De^{m}h=0 &  & \text{in }B_{1}.
\end{eqnarray*}

\textbf{Step 3.} Estimates of $u_{1},u_{2}$ and $h$.

To estimate $u_{1}$, let
\[
u_{11}\equiv\log\ast\left(\sum_{i=1}^{m-1}\de^{i}\left(\sum_{j=m-i}^{m}\na^{j}A\na^{2m-i-j}u\right)+\de K\right)
\]
and
\[
u_{12}=\log\ast\left(Af\right)
\]
such that
\[
u_{1}=u_{11}+u_{12}.
\]
The estimate of $u_{12}$ is simple and standard: since $f\in L^{p}(\R^{2m})$
and $A\in L^{\wq}$, we have $u_{12}\in W^{2m,p}(\R^{2m})$ with $\na^{j}u_{12}=I_{2m-j}(Af)\in L^{p_{j}}(\R^{2m})$
for all $1\le j\le2m$, where $p_{j}^{-1}=p^{-1}-(2m-j)/2m$ and
\begin{equation}
\|\na^{j}u_{12}\|_{L^{p_{j}}(B_{1})}\le\|\na^{j}u_{12}\|_{L^{p_{j}}(\R^{2m})}\le C_{m,j,p}\|f\|_{L^{p}(\R^{2m})}\le C_{m,j,p}\|f\|_{L^{p}(B_{1})}.\label{eq: estimate of u12}
\end{equation}
Here, when $j=2m$, we use $I_{0}$ to denote the singular integral
operators $\na^{2m}\log$ on $L^{p}(\R^{2m})$.

We estimate $u_{11}$ as follows. Note that for each $1\le t\le m$,
we have
\[
\na^{t}u_{11}=I_{2m-t}\left(\sum_{i=1}^{m-1}\de^{i}
\left(\sum_{j=m-i}^{m}\na^{j}A\na^{2m-i-j}u\right)+\de K\right).
\]
Hence
\[
|\na^{t}u_{11}|\lesssim\sum_{i=1}^{m-1}\sum_{j=m-i}^{m}I_{2m-t-i}
\left(|\na^{j}A\na^{2m-i-j}u|\right)+I_{2m-t-1}(|K|).
\]
Since $A,u\in W^{m,2}(\R^{2m})$, we have $\na^{j}A\na^{2m-i-j}u\in L^{\frac{2m}{j},2}\cdot L^{\frac{2m}{2m-i-j},2}\subset L^{\frac{2m}{2m-i},1}\subset L^{\frac{2m}{2m-i},2},$
with
\[
\left\|\na^{j}A\na^{2m-i-j}u\right\|_{L^{\frac{2m}{2m-i},2}}
\le\left\|\na^{j}A\right\|_{L^{\frac{2m}{j},2}}
\left\|\na^{2m-i-j}u\right\|_{L^{\frac{2m}{2m-i-j},2}}
\lesssim\ep\left\|\na^{2m-i-j}u\right\|_{L^{\frac{2m}{2m-i-j},2}}.
\]
Thus, it holds
\[
\left\Vert I_{2m-t-i}\left(|\na^{j}A\na^{2m-i-j}u|\right)\right\Vert _{L^{2m/t,2}(\R^{2m})}\lesssim\ep\left\Vert \na^{2m-i-j}u\right\Vert _{L^{\frac{2m}{2m-i-j},2}}.
\]
Taking summation over $i$ and $j$, and then over $t$, we obtain,
\begin{equation}
\begin{aligned}  &\quad\sum_{t=1}^{m}\left\Vert \sum_{i=1}^{m-1}\sum_{j=m-i}^{m}I_{2m-t-i}\left(|\na^{j}A\na^{2m-i-j}u|\right)\right\Vert _{L^{2m/t,2}(\R^{2m})}\\
 & \lesssim\ep\sum_{j=1}^{m}\|\na^{j}u\|_{L^{\frac{2m}{j},2}(\R^{2m})}  \lesssim\ep\sum_{j=1}^{m}\|\na^{j}u\|_{L^{\frac{2m}{j},2}(B_{1})},
\end{aligned}
\label{eq: estimate of the first term of u11}
\end{equation}
where we used the boundedness of the extensions of $u$ from $B_{1}$
into $\R^{2m}$.

Next consider the first term $K_{1}$ of $K$, i.e.
\[
K_{1}=\sum_{k=1}^{m-1}\sum_{l=0}^{k-1}\De^{l}A\De^{k-l-1}d(V_{k}du).
\]
Repeatedly using Leibniz rule yields
\[
\De^{l}A\De^{k-l-1}d(V_{k}du)=\sum_{i=1}^{2(k-l)}\De^{l}A\na^{2(k-l)-i}V_{k}\na^{i}u.
\]
The estimate of $\left\Vert I_{2m-t-1}\left(\De^{l}A\na^{2(k-l)-i}V_{k}\na^{i}u\right)\right\Vert _{L^{2m/t,2}(\R^{2m})}$
is far more complicated than the previous term. We divide this term
into two cases according to the order of the derivatives of $A$:

\textbf{Case 1:} $2l>m$. Then $2(k-l)<2(m-1-\frac{m}{2})=m-2<m$
and $2(k-l)-1<2k+1-m$. Note that $\De^{l}A\in W^{2l-m,2}$ is of
negative exponent. Since $i\le2(k-l)$, $\na^{2(k-l)-i}V_{k}\in W^{2l-m+i+1,2}$
and $\na^{i}u\in W^{m-i,2}$ are of positive power. Moreover, it is
direct to verify that
\[
\min_{1\le i\le2(k-l)}(2l-m+i+1,m-i)\ge2l-m,
\]
which means, by Proposition \ref{prop: generalized holder ineq}, that the product $\De^{l}A\na^{2(k-l)-i}V_{k}\na^{i}u\in W^{2l-m,2}$.
Using Leibniz rule, we obtain
\begin{equation}
\begin{aligned}\De^{l}A\na^{2(k-l)-i}V_{k}\na^{i}u & =\sum_{|\al|\le2l-m}\pa^{\al}A_{\al}\left(\na^{2(k-l)-i}V_{k}\na^{i}u\right)\\
 & =\sum_{|\al|\le2l-m}\sum_{0\le\be\le\al}\pa^{\be}\left(A_{\al}\pa^{\al-\be}\left(\na^{2(k-l)-i}V_{k}\na^{i}u\right)\right)\\
 & =\sum_{|\al|\le2l-m}\sum_{0\le\be\le\al}\sum_{0\le\ga\le\al-\be}\pa^{\be}\left(A_{\al}\na^{2(k-l)-i+|\ga|}V_{k}\na^{|\al|-|\be|-|\ga|+i}u\right),
\end{aligned}\label{eq: case 1}
\end{equation}
where $A_{\al}\in L^{2}(\R^{2m})$. So $I_{2m-1-t}$ acts on this
summation gives a summation of the form
\[
I_{2m-1-t-|\be|}\left(A_{\al}\na^{2(k-l)-i-|\ga|}V_{k}\na^{|\al|-|\be|-|\ga|+i}u\right).
\]
The integrability of the function inside the potential is as follows: $A_{\al}\in L^{2},$ and
\[
\na^{2(k-l)-i-|\ga|}V_{k}\in W^{2l-m+i+1-\ga,2}\subset L^{\frac{2m}{2m-(2l+i+1)+|\ga|},2},\na^{|\al|-|\be|-|\ga|+i}u\in L^{\frac{2m}{|\al|-|\be|-|\ga|+i},2}
\]
and so
\[
A_{\al}\na^{2(k-l)-i-|\ga|}V_{k}\na^{|\al|-|\be|-|\ga|+i}u\in L^{p,2}
\]
with
\[
\frac{1}{p}=\frac{1}{2}+\frac{2m-(2l+i+1)+|\ga|}{2m}+\frac{|\al|-|\be|-|\ga|+i}{2m}=\frac{3m+|\al|-|\be|-2l-1}{2m}.
\]
Recall that all the functions are extended in such a way that they
are supported in the ball $B_{2}$. Hence we only need $p\ge\frac{2m}{2m-1-t-|\be|+j}$
such that
\[
I_{2m-1-t-|\be|}\left(A_{\al}\na^{2(k-l)-i-|\ga|}V_{k}\na^{|\al|-|\be|-|\ga|+i}u\right)\in L^{2m/t,2}.
\]
This does hold by noticing that
\[
\left(3m+\al-\be-2l-1\right)-(2m-1-\be)=m-2l+\al\le0.
\]
Thus, there holds
\[
\left\Vert I_{2m-1-t-|\be|}\left(A_{\al}\na^{2(k-l)-i-|\ga|}V_{k}\na^{|\al|-|\be|-|\ga|+i}u\right)\right\Vert _{L^{2m,2}}\lesssim\ep\sum_{i=1}^{m}\|\na^{i}u\|_{L^{2m/i,2}(B_{1})}.
\]
In conclusion, in \textbf{Case 1} we get
\begin{equation}
\left\Vert I_{2m-t-1}\left(\De^{l}A\na^{2(k-l)-i}V_{k}\na^{i}u\right)\right\Vert _{L^{2m/t,2}(\R^{2m})}\lesssim\ep\sum_{i=1}^{m}\|\na^{i}u\|_{L^{2m/i,2}(B_{1})}.\label{eq: estamate for Case 1}
\end{equation}

\textbf{Case 2}: $2l\le m$. In this case, we have two sub-cases:

\textbf{Subcase 2.1:} $2l\le m$ and $2(k-l)-i>2k+1-m$. That is,
$\na^{2(k-l)-i}V_{k}\in W^{2l+i+1-m,2}$ is of negative exponent.
Then, $i<m-2l-1<m$. In this case
\[
\min_{1\le i\le2(k-l)}(m-2l,m-i)\ge|2l+i+1-m|=m-(2l+i+1).
\]
So $\De^{l}A\na^{2(k-l)-i}V_{k}\na^{i}u$ makes sense, and
\begin{equation}
\begin{aligned}\De^{l}A\na^{2(k-l)-i}V_{k}\na^{i}u & =\sum_{|\al|\le m-(2l+i+1)}\De^{l}A\pa^{\al}V_{k,\al}\na^{i}u\\
 & =\sum_{\al}\sum_{\be}\pa^{\be}\left(V_{k,\al}\pa^{\al-\be}(\De^{l}A\na^{i}u)\right)\\
 & =\sum_{\al}\sum_{\be}\sum_{\ga}\pa^{\be}\left(V_{k,\al}
 \na^{2l+|\ga|}A\na^{|\al|-|\be|-|\ga|+i}u\right).
\end{aligned}
  \label{eq: case 2.1}
\end{equation}
Then
\[
V_{k,\al}\na^{2l+|\ga|}A\na^{|\al|-|\be|-|\ga|+i}u\in L^{p,2}
\]
with
\[
\frac{1}{p}=\frac{1}{2}+\frac{2l+|\ga|}{2m}+\frac{|\al|-|\be|-|\ga|+i}{2m}=\frac{m+2l+|\al|-|\be|+i}{2m}.
\]
Thus $p\ge\frac{2m}{2m-1-t-\be+t}$ is equivalent to
\[
m+2l+|\al|-|\be|+i\le2m-1-|\be|,
\]
which is equivalent to
\[
\al\le m-2l-i-1.
\]
This holds by our assumption. So in \textbf{Subcase 2.1} we also have
the estimate (\ref{eq: estamate for Case 1}).

\textbf{Subcase 2.2:} $2l\le m$ and $2(k-l)-i\le2k+1-m$. That is,
$\De^{l}A$ and $\na^{2(k-l)-i}V_{k}\in W^{2l+i+1-m,2}$ are  usual
Sobolev functions with positive exponent. In this case we have two more
cases:

(1) $i\le m$. Then
\[
\De^{l}A\na^{2(k-l)-i}V_{k}\na^{i}u\in L^{p,2}
\]
with
\[
\frac{1}{p}=\frac{2l}{2m}+\frac{2m-(2l+i+1)}{2m}+\frac{i}{2m}=\frac{2m-1}{2m}=\frac{2m-1-t+t}{2m}.
\]
Then $I_{2m-1-t}(L^{\frac{2m}{2m-1},q})\subset L^{2m,2}$ holds.

(2) $i>m$. Then $\na^{i}u\in W^{m-i,2}$ is of negative exponent.
In this case
\[
\min(m-2l,2l+i+1-m)\ge i-m.
\]
So
\begin{equation}
\begin{aligned}\De^{l}A\na^{2(k-l)-i}V_{k}\na^{i}u & =\sum_{\al\le i-m}\De^{l}A\na^{2(k-l)-i}V_{k}\pa^{\al}u_{\al}\\
 & =\sum_{\al\le i-m}\sum_{\be}\pa^{\be}\left(u_{\al}\pa^{\al-\be}\left(\De^{l}A\na^{2(k-l)-i}V_{k}\right)\right)\\
 & =\sum_{\al\le i-m}\sum_{\be}\sum_{\ga}\pa^{\be}\left(u_{\al}\na^{2l+|\ga|}A\na^{2(k-l)-i+|\al|-|\be|-|\ga|}V_{k}\right).
\end{aligned}
  \label{eq: case 2.2}
\end{equation}
We have
\[
u_{\al}\na^{2l+|\ga|}A\na^{2(k-l)-i+|\al|-|\be|-|\ga|}V_{k}\in L^{p,2}
\]
with
\[
\frac{1}{p}=\frac{1}{2}+\frac{2l+|\ga|}{2m}+\frac{2m-1-2l-i+(|\al|-|\be|-|\ga|)}{2m}=\frac{3m-1+i+|\al|-|\be|}{2m}.
\]
Thus $p\ge\frac{2m}{2m-2-|\be|+1}$ is equivalent to
\[
3m-1-i+|\al|-|\be|\le2m-1-|\be|,
\]
is equivalent to
\[
\al\le i-m,
\]
which holds. Therefore, in \textbf{Subcase 2.2}  the
estimate (\ref{eq: estamate for Case 1}) holds as well.

In conclusion, we obtain the estimate
\begin{equation}
\sum_{l=1}^{m}\left\Vert I_{2m-1-t}\left(\sum_{k=1}^{m-1}\sum_{l=0}^{k-1}\De^{l}A\na^{2(k-l)-i}V_{k}\na^{i}u\right)\right\Vert _{L^{2m/t,2}(\R^{2m})}\lesssim\ep\sum_{i=1}^{m}\|\na^{i}u\|_{L^{2m/i,2}(B_{1})}.\label{eq: estimate of K1}
\end{equation}
Arguing exactly as in the above, we can obtain the same estimate for
the remaining terms of $K$. Therefore, combining (\ref{eq: estimate of the first term of u11})
and (\ref{eq: estimate of K1}) together, we deduce
\begin{equation}
\sum_{j=1}^{m}\|\na^{j}u_{11}\|_{L^{2m/j,2}(\R^{2m})}\lesssim\ep\sum_{i=1}^{m}\|\na^{i}u\|_{L^{2m/i,2}(B_{1})}.\label{eq: estimate of u11}
\end{equation}
Using (\ref{eq: estimate of u11}) and (\ref{eq: estimate of u12}),
we obtain
\begin{equation}
\sum_{j=1}^{m}\|\na^{j}u_{1}\|_{L^{2m/j,2}(\R^{2m})}\lesssim\ep\sum_{i=1}^{m}\|\na^{i}u\|_{L^{2m/i,2}(B_{1})}+\|f\|_{L^{p}(B_{1})}.\label{eq: estimates of u1}
\end{equation}

Next we estimate the function $u_{2}$ and $h$. By the definition
(\ref{eq: def of u2}) of $u_{2}$, we use the same method as that
of (\ref{eq: estimate of the first term of u11}) to obtain
\begin{equation}
\sum_{j=1}^{m}\|\na^{j}u_{2}\|_{L^{2m/j,2}(\R^{2m})}\lesssim\ep\sum_{i=1}^{m}\|\na^{i}u\|_{L^{2m/i,2}(B_{1})}.\label{eq: estimate of u2}
\end{equation}
As to the polyharmonic function $h$, we can apply \cite[Lemma 6.2]{Gastel-Scheven-2009CAG} to find that,
for any $0<r<1$,
\begin{equation}
\sum_{j=1}^{m}\|\na^{j}h\|_{L^{2m/j,2}(B_{r})}\lesssim r\sum_{i=1}^{m}\|\na^{i}h\|_{L^{2m/i,2}(B_{1})}.\label{eq: decay of h}
\end{equation}

\textbf{Step 4.} Conclusion.

Now we prove the decay estimate \eqref{eq: decay estimates of higher order system} as follows. Let $0<\tau<1$
be determined later. We have, for every $1\le j\le m$
\[
\begin{aligned}\|\na^{j-1}\left(A^{-1}du_{11}\right)\|_{L^{2m/j,2}(B_{\tau})} & \lesssim\sum_{i=1}^{j}\|\na^{i}u_{11}\|_{L^{\frac{2m}{i},2}(B_{\tau})}\\
 & \begin{aligned} & \lesssim\sum_{i=1}^{j}\|\na^{i}u_{11}\|_{L^{\frac{2m}{i},2}(B_{1})}
 \lesssim\ep\sum_{i=1}^{m}\|\na^{i}u\|_{L^{2m/i,2}(B_{1})};\end{aligned}
\end{aligned}
\]
\[
\begin{aligned}\|\na^{j-1}\left(A^{-1}du_{12}\right)\|_{L^{2m/j,2}(B_{\tau})} & \lesssim\sum_{i=1}^{j}\|\na^{i}u_{12}\|_{L^{\frac{2m}{i},2}(B_{\tau})}\\
 & \lesssim\tau^{2m(1-\frac{1}{p})}\sum_{i=1}^{j}\|\na^{i}u_{12}\|_{L^{p_{j}}(B_{\tau})}
 \le\tau^{2m(1-\frac{1}{p})}\|f\|_{L^{p}(B_{1})};
\end{aligned}
\]
\[
\|\na^{j-1}\left(A\ast du_{2}\right)\|_{L^{2m/j,2}(B_{\tau})}\lesssim\sum_{i=1}^{j}\|\na^{i}u_{2}\|_{L^{\frac{2m}{i},2}(B_{\tau})}
\lesssim\ep\sum_{i=1}^{m}\|\na^{i}u\|_{L^{2m/i,2}(B_{1})};
\]
\[
\|\na^{j-1}(A^{-1}h)\|_{L^{2m/j,2}(B_{\tau})}\lesssim\sum_{i=1}^{j}\|\na^{i}h\|_{L^{\frac{2m}{i},2}(B_{\tau})}
\lesssim\tau\sum_{i=1}^{m}\|\na^{i}h\|_{L^{2m/i,2}(B_{1})}.
\]
Therefore, from the above four estimates we derive, for each $1\le j\le m$,
\[
\begin{aligned}\|\na^{j}u\|_{L^{\frac{2m}{j},2}(B_{\tau})} & \le\|\na^{j-1}(A^{-1}h)\|_{L^{\frac{2m}{j},2}(B_{\tau})}
+\|\na^{j-1}\left(A^{-1}du_{1}\right)\|_{L^{\frac{2m}{j},2}(B_{\tau})}\\
&\quad+\|\na^{j-1}\left(A\ast du_{2}\right)\|_{L^{\frac{2m}{j},2}(B_{\tau})}\\
 & \lesssim\tau\sum_{i=1}^{m}\|\na^{i}h\|_{L^{2m/i,2}(B_{1})}+\ep\sum_{i=1}^{m}\|\na^{i}u\|_{L^{2m/i,2}(B_{1})}
 +\tau^{2m(1-\frac{1}{p})}\|f\|_{L^{p}(B_{1})}\\
 & \lesssim\tau\sum_{i=1}^{m}\left(\|\na^{i}u\|_{L^{2m/i,2}(B_{1})}+\|\na^{i}u_{11}\|_{L^{2m/i,2}(B_{1})}
 +\sum_{i=1}^{m}\|\na^{i}u_{12}\|_{L^{2m/i,2}(B_{1})}\right)\\
 & \quad+\ep\sum_{i=1}^{m}\|\na^{i}u\|_{L^{2m/i,2}(B_{1})}+\tau^{2m(1-\frac{1}{p})}\|f\|_{L^{p}(B_{1})}\\
 & \le C(\tau+\ep)\sum_{i=1}^{m}\|\na^{i}u\|_{L^{2m/i,2}(B_{1})}
 +C\left(\tau+\tau^{2m(1-\frac{1}{p})}\right)\|f\|_{L^{p}(B_{1})}
\end{aligned}
\]
for some $C>0$ independent of $\tau$ and $\ep$.

Now set
$
\al_{0}=2m(1-1/p)
$ and
 $\be=(\al_{0}+1)/2$. First choose $\tau$ such that $2C\tau\le\tau^{\be}$,
and then choose $\ep<\tau$. Taking summation over $j$, we obtain
\[
\sum_{i=1}^{m}\|\na^{i}u\|_{L^{2m/i,2}(B_{\tau})}\le\tau^{\be}\sum_{i=1}^{m}\|\na^{i}u\|_{L^{2m/i,2}(B_{1})}
+C\tau^{\al_{0}}\|f\|_{L^{p}(B_{1})}
\]
for some $C>0$ independent of $\tau$. This proves \eqref{eq: decay estimates of higher order system}.

Finally, apply an iteration argument as follows. Write $\Phi(u,B_{\tau})=\sum_{i=1}^{m}\|\na^{i}u\|_{L^{\frac{2m}{i},2}(B_{\tau})}$
and $F=\|f\|_{L^{p}(B_{1})}$. Then scaling invariance implies
\[
\Phi(u_{\tau},B_{1})=\Phi(u,B_{\tau})\qquad\text{and}
\qquad\|f_{\tau}\|_{L^{p}(B_{1})}=\tau^{\al_{0}}\|f\|_{L^{p}(B_{\tau})}.
\]
Thus, for any $k\in\N$, there holds
\[
\begin{aligned}\Phi(u,B_{\tau^{k}})=\Phi(u_{\tau^{k-1}},B_{\tau}) & \le\tau^{\be}\Phi(u_{\tau^{k-1}},B_{1})+C\tau^{\al_{0}}\|f_{\tau^{k-1}}\|_{L^{p}(B_{1})}\\
 & =\tau^{\be}\Phi(u,B_{\tau^{k-1}})+C\tau^{k\al_{0}}\|f\|_{L^{p}(B_{\tau^{k-1}})}\\
 & \le\tau^{\be}\Phi(u,B_{\tau^{k-1}})+C\|f\|_{L^{p}(B_{1})}\tau^{k\al_{0}}.
\end{aligned}
\]
Therefore, using an iteration arguments we obtain
\[
\Phi(u,B_{\tau^{k}})\le\tau^{k\be}\Phi(u,B_{1})+C\|f\|_{L^{p}(B_{1})}\tau^{k\al_{0}}
\sum_{i=0}^{k-1}\tau^{i(\be-\al_{0})}.
\]
Since $\be>\al_{0}$, we obtain, for any $k\ge1$,
\[
\Phi(u,B_{\tau^{k}})\le C\tau^{k\al_{0}}\left(\Phi(u,B_{1})+\|f\|_{L^{p}(B_{1})}\right).
\]
 The proof of Lemma \ref{lem: decay estimate lemma 2}  is complete taking into account of the monotonicity of
$r\mapsto\Phi(u,B_{r})$.

\section{Higher order regularity}\label{sec:higher order Sobolev regularity}

In this section, we will prove an almost optimal  higher order regularity result. The main tool is   Lemma \ref{lemma:improved Riesz potential}.  Throughout this section, we set $$1<p<\frac{2m}{2m-1}, \quad
\al_0=2m(1-1/p) \quad \text{ and } \quad M\equiv \|u\|_{W^{m,2}(B_1)}+\|f\|_{L^p(B_1)}.$$

\subsection{$W^{m,q}$-estimate with any $q<\frac{2p}{2-p}$}

\begin{proposition}\label{prop:m order Sobolev regularity}
Suppose $f\in L^p(B_{1})$ for some $p\in (1,{2m}/{2m-1})$ and $u\in W^{m,2}(B_{1},\R^n)$ is a  solution of the inhomogeneous system \eqref{eq:nonhomo Longue-Gastel system}.  Then $u\in W^{m,q}_{\loc}(B_{1})$  for any $q<\frac{2p}{2-p}$.
\end{proposition}

\begin{proof}
As in the previous section, it suffices to assume that $\ta_{B_1}<\ep$ for some sufficiently small $\ep$, such that the conservation law \eqref{eq:conservation law of D-G} holds for some functions $A,B$. And then extend all the related functions from $B_{1/2}$ to the whole space with compact support in $B_2$ in a bounded way. Then apply the Hodge decomposition for $Adu$ as follows.

Set
\[
v=\log\ast\left(\sum_{i=1}^{m-1}\de^{i}\left(\sum_{j=m-i}^{m}\na^{j}A\na^{2m-i-j}u\right)+\de K+Af\right),
\]
\[
g=\log\ast\De^{m-1}(dA\wedge du)
\]
and
\[
h=Adu-dv-\ast dg
\]
such that $h$ is $m$-polyharmonic in $B_{1}$.

To begin with, recall that we have proved in Theorem  \ref{thm:optimal Holder exponent for inho DG} that
\[
\na^{i}u\in M^{\frac{2m}{i},\frac{2m}{i}\al_0}(\R^{2m})\cap L^{\frac{2m}{i}}(\R^{2m}),\qquad i=1,\ldots,m.
\]



\textbf{Estimate of $g$.}
By H\"older's inequality, for every $1\le i\le m$,
\[
\na^{i}A\na^{m+1-i}u\in L^{\frac{2m}{m+1}}\cap M^{1,m-1+\al_0}(\R^{2m})=L^{\frac{2m}{m+1}}\cap M^{1,2m-(m+1-\al_0)}(\R^{2m}).
\]
Thus $$\na^{m}g\approx I_{1}\left(\sum_{i=1}^{m}\na^{i}A\na^{m+1-i}u\right)\in L^{q_{0}}(\R^{2m})$$
with
\begin{equation}\label{eq: q_0}
q_{0}=\frac{2m}{m+1}\frac{m+1-\al_0}{m-\al_0}>2
\end{equation}
since $0<\al_0<m$.

To estimate $v$, we introduce the following notations to decompose $\na^m v$:
\begin{itemize}
\item $v_1:=I_{m}\ast\sum_{i=1}^{m-1}\de^{i}\left(\sum_{j=m-i}^{m}\na^{j}A\na^{2m-i-j}u\right)$;
\item $v_2:=I_m\ast \delta(K)=I_{m-1}(K)$;
\item $v_3:=I_m\ast (Af)$.
\end{itemize}

\textbf{Estimate of $v_1$.}
Note that
\[
v_{1}= I_{m}\sum_{i=1}^{m-1}\de^{i}\left(\sum_{j=m-i}^{m}\na^{j}A\na^{2m-i-j}u\right)\approx\sum_{i=1}^{m-1}\sum_{j=m-i}^{m}I_{m-i}\left(\na^{j}A\na^{2m-i-j}u\right)
\]
and that
\[
\na^{j}A\na^{2m-i-j}u\in L^{\frac{2m}{2m-i}}\cap M^{1,i+\al_0}=L^{\frac{2m-i}{2m}}\cap M^{1,2m-(2m-i-\al_0)}(\R^{2m}).
\]
Hence $I_{m-i}\left(\na^{j}A\na^{2m-i-j}u\right)\in L^{q_{i}}$ with
\[
\frac{1}{q_{i}}=\frac{2m-i}{2m}\left(1-\frac{m-i}{2m-i-\al_0}\right)
=\frac{(2m-i)(m-\al_0)}{2m(2m-i-\al_0)}<\frac{1}{2}
\]
for all $1\le i\le m-1$. Moreover, we have $q_{1}>q_{2}>\cdots>q_{m-1}=q_0>2$.
Thus we obtain
\[
v_{1}\in L^{q_{m-1}}=L^{q_{0}}(\R^{2m}),
\]
where $q_{0}$ is the number defined as in (\ref{eq: q_0}).

\textbf{Estimate of $v_2$.}
Next we shall establish similar estimates for the terms involving $K$. Let $K_{1}$ be the first term  of $K$, i.e.
\[
K_{1}=\sum_{k=1}^{m-1}\sum_{l=0}^{k-1}\De^{l}A\De^{k-l-1}d(V_{k}du)
=\sum_{k=1}^{m-1}\sum_{l=0}^{k-1}\sum_{i=1}^{2(k-l)}\De^{l}A\na^{2(k-l)-i}V_{k}\na^{i}u.
\]
by repeatedly using Leibniz rules.
As in the proof of Lemma \ref{lem: decay estimate lemma 2},  we divide this term into two cases.

\textbf{Case 1:} $2l>m$. Then by \eqref{eq: case 1} there holds
\[
\De^{l}A\na^{2(k-l)-i}V_{k}\na^{i}u =\sum_{|\al|\le2l-m}\sum_{0\le\be\le\al}\sum_{0\le\ga\le\al-\be}\pa^{\be}
\left(A_{\al}\na^{2(k-l)-i+|\ga|}V_{k}\na^{|\al|-|\be|-|\ga|+i}u\right),
\]
for some $A_{\al}\in L^{2}(\R^{2m})$. So $I_{m-1-|\beta|}$ acts on this
term giving a summation of the type
\[
I_{m-1-|\be|}\left(A_{\al}\na^{2(k-l)-i-|\ga|}V_{k}\na^{|\al|-|\be|-|\ga|+i}u\right).
\]
The integrability of the function inside the potential is as follows:
\[
A_{\al}\in L^{2},\na^{2(k-l)-i-|\ga|}V_{k}\in W^{2l-m+i+1-\ga}\subset L^{\frac{2m}{2m-(2l+i+1)+|\ga|}},\na^{|\al|-|\be|-|\ga|+i}u\in L^{\frac{2m}{|\al|-|\be|-|\ga|+i}}
\]
and so
\[
T:=A_{\al}\na^{2(k-l)-i-|\ga|}V_{k}\na^{|\al|-|\be|-|\ga|+i}u\in L^{p}
\]
with
\[
\frac{1}{p}=\frac{1}{2}+\frac{2m-(2l+i+1)+|\ga|}{2m}+\frac{|\al|-|\be|-|\ga|+i}{2m}=\frac{3m+|\al|-|\be|-2l-1}{2m}.
\]
Note that since $|\alpha|-2l\leq -m$,
$$3m+|\alpha|-|\beta|-2l-1\leq 2m-1-|\beta|$$
and so $T\subset L^{p}\subset L^{\frac{2m}{2m-1-|\beta|}}$.

By H\"older's inequality and the decay estimate in Theorem \ref{thm:optimal Holder exponent for inho DG}, we can  verify that
\[
\sup_{x\in B_{1/2},0<r<1/2}r^{-q\al_0}\int_{B_{r}(x)}|T|^{q}\le CM,
\]
where $q=q_\beta=\frac{2m}{2m-1-|\beta|}$. By H\"older's inequality, this implies  that $K_1\in M^{1,\alpha_\beta+1}(B_{1/2})$, that is,
\[
\sup_{x\in B_{1/2},0<r<1/2}r^{-(1+\al_\beta)}\int_{B_{r}(x)}|T|\le CM,
\]
where $1+\alpha_\beta=\alpha_0+2m(1-1/q_\beta)=\alpha_0+1+|\beta|$. Hence
$$I_{m-1-|\beta|}\left(T \right)\in L^{q_\beta\frac{2m-1-|\beta|-\alpha_0}{2m-1-|\beta|-\alpha_0-(m-1-|\beta|)}}
=L^{\frac{2m}{2m-1-|\beta|}\frac{2m-1-|\beta|-\alpha_0}{m-\alpha_0}}.$$
Note that as a function of $\beta$, $\frac{2m}{2m-1-|\beta|}\frac{2m-1-|\beta|-\alpha_0}{m-\alpha_0}$ is decreasing and
$$|\beta|\leq |\al|\le  2l-m\leq 2(m-2)-m=m-4\leq m-2.$$  This implies
\[\frac{2m}{2m-1-|\beta|}\frac{2m-1-|\beta|-\alpha_0}{m-\alpha_0} \ge q_0\] for all $|\be|\le |\al| \le  2l-m$.
So
$$I_{m-1}(K_1)\approx I_{m-1-|\beta|}\left(T \right)\in L^{q_0}.$$



\textbf{Case 2}: $2l\le m$. In this case, we have two subcases:

\textbf{Subcase 2.1:} $2l\le m$ and $2(k-l)-i>2k+1-m$. Then, by \eqref{eq: case 2.1} there holds
\[
\De^{l}A\na^{2(k-l)-i}V_{k}\na^{i}u =\sum_{\al}\sum_{\be}\sum_{\ga}\pa^{\be}
\left(V_{k,\al}\na^{2l+|\ga|}A\na^{|\al|-|\be|-|\ga|+i}u\right).
\]
Thus,
\[
T:=V_{k,\al}\na^{2l+|\ga|}A\na^{|\al|-|\be|-|\ga|+i}u\in L^{p}
\]
with
\[
\frac{1}{p}=\frac{1}{2}+\frac{2l+|\ga|}{2m}+\frac{|\al|-|\be|-|\ga|+i}{2m}=\frac{m+2l+|\al|-|\be|+i}{2m}.
\]
Since $|\alpha|\leq m-(2l+i+1)$,
$$m+2l+|\alpha|-|\beta|-1\leq 2m-1-|\beta|$$
and we are in the same situation as in \textbf{Case 1}, that is,
$$T\in L^{\frac{2m}{2m-1-|\beta|}}\cap M^{1,\alpha_0+1+|\beta|}.$$
So $$I_{m-1-|\beta|}\left(T \right)\in L^{\frac{2m}{2m-1-|\beta|}\frac{2m-1-|\beta|-\alpha_0}{m-\alpha_0}}.$$
Note that
$$|\beta|\leq m-(2l+i+1)\leq m-i-1\leq m-2.$$
 Thus we have
$$I_{m-1}(K_1)\approx I_{m-1-|\beta|}(T)\in L^{q_0}.$$

\textbf{Subcase 2.2:} $2l\le m$ and $2(k-l)-i\le2k+1-m$. That is, $\De^{l}A$
and $\na^{2(k-l)-i}V_{k}\in W^{2l+i+1-m,2}$ are the usual Sobolev
functions of positive exponents. In this case we have two more cases:
(1) $i\le m$ and (2) $i>m$.

(1). When $i\le m$, we have
\[
\De^{l}A\na^{2(k-l)-i}V_{k}\na^{i}u\in L^{p,2}
\]
with
\[
\frac{1}{p}=\frac{2l}{2m}+\frac{2m-(2l+i+1)}{2m}+\frac{i}{2m}=\frac{2m-1}{2m}.
\]
Then $I_{m-1}(\De^{l}A\na^{2(k-l)-i}V_{k}\na^{i}u)\in L^{\frac{2m}{2m-1}\frac{2m-1-\alpha_0}{2m-1-\alpha_0-(m-1)}}\subset L^{q_0}$ holds.

(2). When $i>m$, by \eqref{eq: case 2.2} we have
\[
\De^{l}A\na^{2(k-l)-i}V_{k}\na^{i}u  =\sum_{\al\le i-m}\sum_{\be}\sum_{\ga}\pa^{\be}
\left(u_{\al}\na^{2l+|\ga|}A\na^{2(k-l)-i+|\al|-|\be|-|\ga|}V_{k}\right).
\]
So
\[
T:=u_{\al}\na^{2l+|\ga|}A\na^{2(k-l)-i+|\al|-|\be|-|\ga|}V_{k}\in L^{p,2}
\]
with
\[
\frac{1}{p}=\frac{1}{2}+\frac{2l+|\ga|}{2m}+\frac{2m-1-2l-i+(|\al|-|\be|-|\ga|)}{2m}=\frac{3m-1-i+|\al|-|\be|}{2m}.
\]
Since $|\alpha|\leq i-m$,
\[
3m-1-i+|\al|-|\beta|\le 2m-1-|\beta|
\]
and so
$$T\in L^{\frac{2m}{2m-1-|\beta|}}\cap M^{1,\alpha_0+1+|\beta|}$$
and
$$I_{m-1-|\beta|}\left(T \right)\in L^{\frac{2m}{2m-1-|\beta|}\frac{2m-1-|\beta|-\alpha_0}{m-\alpha_0}}.$$
Note that
\[
\beta\leq i-m\leq 2k-m\leq 2(m-1)-m=m-2,
\]
which again is similar to the \textbf{Case 1}. Therefore, in \textbf{Subcase 2.2} we conclude that
$$I_{m-1}(K_1)\approx I_{m-1-|\beta|}(T)\in L^{q_0}.$$

Combining all the above estimate together, we  conclude that
$I_{m-1}(K_1)\in L^{q_0}.$
One can estimate similarly for the other terms in $K$ to arrive finally at
$$v_2\in L^{q_0}.$$

\textbf{Estimate of $v_3$.}
It follows from standard elliptic regularity theory that
\[
v_3\in W^{m,p}(\R^{2m})\subset L^{\frac{2p}{2-p}}(\R^{2m}).
\]
In conclusion, we obtain \[
v\in   L^{q_0}(\R^{2m})\] for some $q_0>2$.

 Since the polyharmonic function $h$ is smooth in $B_1$,  we deduce
\[u\in W_{\loc}^{m,q_{0}}(B_{1}).\]
Note that  $p>1$ and $\al_0>0$ implies $2<q_{0}<\frac{2p}{2-p}$.  Thus we have improved the regularity of $u$ from $W^{m,2}$ to $W^{m, q_0}$.

\textbf{Iteration.}
Next we use a bootstrapping argument to continue improving the regularity of $u$.
We claim that
\begin{eqnarray}\label{eq: bootstrapping W2,q}
\begin{aligned}
u\in W^{m,q}_{loc}\quad \text{with }q<\frac{2p}{2-p} \Longrightarrow u\in W^{m,\frac{2mq_0}{q_0+2m}\frac{m+1-\alpha}{m-\alpha}}_{loc}.
\end{aligned}
\end{eqnarray}
This is true because if $u\in W^{m,q}_{\loc}$ with $q<\frac{2p}{2-p}$, then
\begin{equation}\label{eq:claim for iteration}
v_0,v_1,v_2\in L^{\frac{2mq_0}{q_0+2m}\frac{m+1-\alpha_0}{m-\alpha_0}}.
\end{equation}
Indeed, consider first $v_0$ and we know
$$\nabla^{m+1-i}u\in W^{i-1,q_0}\subset L^{\frac{2mq_0}{2m-(i-1)q_0}}$$
and so $T=\sum_i\nabla^iA\nabla^{m+1-i}u\in L^{\frac{2mq_0}{2m+q_0}}\cap M^{1,m-1+\alpha_0}$. Lemma \ref{lemma:improved Riesz potential} implies
$$v_0\approx I_{1}(T)\in L^{\frac{2mq_0}{2m+q_0}\frac{m+1-\alpha_0}{m-\alpha_0}}.$$
For the term $v_1$, we have
$$\nabla^{2m-i-j}u\in W^{m-(2m-i-j),q_0}\subset L^{\frac{2mq_0}{2m-(i+j-m)q_0}}$$
and so $T=\nabla^jA\nabla^{2m-i-j}u\in L^{\tilde{q}_{i}}$ with $\tilde{q}_{i}=\frac{2mq_0}{(m-i)q_0+2m}$.
Since $T\in M^{1,1+\al_0}(B_{\frac{1}{2}})$, Lemma \ref{lemma:improved Riesz potential} implies that
\[
v_1\approx I_{m-i}(T)\in L^{\tilde{q}_{i}\frac{2m-1-\al_0}{m-\al_0}}
=L^{\frac{2mq_0}{(m-i)q_0+2m}\frac{2m-i-\alpha_0}{m-\alpha_0}}.
\]
Note that as a function of $i$, $\frac{2mq_0}{(m-i)q_0+2m}\frac{2m-i-\alpha_0}{m-\alpha_0}$ is decreasing, and so it attains the minimum $\frac{2mq_0}{q_0+2m}\frac{m+1-\alpha}{m-\alpha}$   when $i=m-1$, and thus
$$v_1 \subset L^{\frac{2mq_0}{q_0+2m}\frac{m+1-\alpha_0}{m-\alpha_0}}.$$
Finally, we shall consider $I_{m-1}(K_1)$ case by case as in the previous study. Corresponding to \textbf{Case 1} above, we have
$$T=A_{\al}\na^{2(k-l)-i-|\ga|}V_{k}\na^{|\al|-|\be|-|\ga|+i}u\in L^p$$
with $p=\frac{2mq_0}{2m+(2m-2l+|\alpha|-|\beta|-1)q_0}$. Since in this case, $2m-2l+|\alpha|-1\leq 2m-m-1=m-1$, we further conclude that
$$T\in L^p\subset L^{\frac{2mq_0}{2m+(m-1-|\beta|)q_0}}.$$
Since $T\in M^{1,\alpha_0+1+|\beta|}$ as well, Lemma \ref{lemma:improved Riesz potential} gives
$$I_{m-1-|\beta|}(T)\in L^{\frac{2mq_0}{2m+(m-1-|\beta|)q_0}\frac{2m-\alpha_0-1-|\beta|}{m-\alpha_0}}.$$
Since as a function of $\beta$, $\frac{2mq_0}{2m+(m-1-|\beta|)q_0}\frac{2m-\alpha_0-1-|\beta|}{m-\alpha_0}$ is decreasing, it attains the minimum $\frac{2mq_0}{2m+q_0}\frac{m+1-\alpha_0}{m-\alpha_0}$ when $|\beta|=m-2$. This gives
$$I_{m-1}(K_1)\in L^{\frac{2mq_0}{2m+q_0}\frac{m+1-\alpha_0}{m-\alpha_0}}.$$ Similarly, one can prove the other cases as well. For instance, in \textbf{Subcase 2.1}, we know
$$\nabla^{|\alpha|-|\beta|-|\gamma|+i}u\in L^{\frac{2mq_0}{2m-(m-|\alpha|+|\beta|+|\gamma|-i)q_0}}$$
and thus
$$T=V_{k,\al}\na^{2l+|\ga|}A\na^{|\al|-|\be|-|\ga|+i}u\in L^{p}$$
with $p=\frac{2mq_0}{2m+(2l+|\alpha|-|\beta|+i)q_0}$. Since
$$2l+|\alpha|-|\beta|+i\leq 2l+m-(2l+i+1)-|\beta|+i=m-1-|\beta|,$$
$T\in L^{\frac{2mq_0}{2m+(m-1-|\beta|)q_0}}\cap M^{1+\alpha_0+|\beta|}$ and so similar to the reason for \textbf{Case 1}, we conclude
$$I_{m-1}(K)\in L^{\frac{2mq_0}{2m+q_0}\frac{m+1-\alpha_0}{m-\alpha_0}}.$$
All together leads to \eqref{eq:claim for iteration}.

Since $v_3\in L^{\frac{2p}{2-p}}$, the claim \eqref{eq:claim for iteration} implies $\nabla^m u\in L^{\frac{2mq_0}{q_0+2m}\frac{m+1-\alpha_0}{m-\alpha_0}}_{\loc}(B_1)$. That is,  $$u\in W^{m,\frac{2mq_0}{q_0+2m}\frac{m+1-\alpha_0}{m-\alpha_0}}.$$

Finally, noticing that
\[\frac{2mq_0}{q_0+2m}\frac{m+1-\alpha_0}{m-\alpha_0}<\frac{2p}{2-p}\Longleftrightarrow q_0<\frac{2p}{2-p}\]
and that  $$q_0\nearrow \frac{2p}{2-p} \quad  \Rightarrow \quad p\frac{2mq_0}{q_0+2m}\frac{m+1-\alpha_0}{m-\alpha_0}\nearrow \frac{2p}{2-p}.$$
Thus, by iterating the bootstrapping claim \eqref{eq: bootstrapping W2,q}, we eventually find that
$$u\in W^{m,q}_{\loc}\quad \text{for all }q<\frac{2p}{2-p}.$$
The proof of Proposition \ref{prop:m order Sobolev regularity} is complete.
\end{proof}


\subsection{$W^{m+1,q}$-estimate with any $q<\frac{2mp}{2m-(m-1)p}$}

\begin{proposition}\label{prop:m+1 order Sobolev regularity}
Let $u\in W^{m,2}(B_{1},\R^n)$ be a weak solution of the inhomogenuous system \eqref{eq:Longue-Gastel system} with $f\in L^p(B_{1})$ for $p\in (1,\frac{2m}{2m-1})$.  Then $u\in W^{m+1,q}_{\loc}(B_{1})$  for any $q<\frac{2mp}{2m-(m-1)p}$.
\end{proposition}
\begin{proof}
Note that $u\in W^{m,q}$ implies $\nabla^i u\in L^{\frac{2mq}{2m-(m-i)q}}$ and so the proof of Proposition \ref{prop:m order Sobolev regularity} implies
$$\na^m g,v_1,v_2,v_3\in W^{1,\frac{2mq}{2m+q}},$$
or equivalently $u\in W^{m+1,\frac{2mq}{2m+q}}$. Note that  $\frac{2mq}{2m+q}\nearrow \frac{2mp}{2m-(m-1)p}$ when $q\nearrow \frac{2p}{2-p}$. The proof is complete.
\end{proof}

\begin{remark}\label{rmk:on best m+1 regularity}
	Once we prove that $u$ belongs to $W^{m,\frac{2p}{2-p}}_{\loc}$, then the above proof implies that $u\in W^{m+1,\frac{2mp}{2m-(m-1)p}}_{\loc}$.
\end{remark}

The following example  shows that the $(m+1)$-th order (Sobolev) regularity is the best one can expect for the general system \eqref{eq:Longue-Gastel system}.

\begin{example}[Solutions without $W^{m+2,p}$-regularity]\label{example:no Wm+2 estimate}
Let $g\colon \R\to \R$ be a continuous function with the following properties:
\begin{itemize}
\item $g\in W^{m+1,2}\big((-2,2)\big)$ but $g\not\in W^{m+2,1}\big((-1,1)\big)$;
\item $g\geq 1$ on $(-1,1)$.
\end{itemize}

Consider the map $u\colon B_1\to \R$, $B_1\subset \R^{2m}$, defined by
$$u(x)=x_1g(x_2).$$
Set
$$V_1(x)=x_1\frac{g''(x_2)}{g(x_2)}\quad \text{and}\quad V(x)=(V_1(x),0,\cdots,0).$$
It is straightforward to verify that $V\in W^{m-1,2}\big(B_1\big)$ and
\begin{equation*}
\Delta^m u=\Delta^{m-1}\big(V\cdot \nabla u\big)\quad \text{in } B_1.
\end{equation*}
However, the regularity of $g$ implies that $u\notin W^{m+2,1}(B_1)$.

\end{example}

\section{Optimal local estimates}\label{sec:optimal global estimates}

In this section, we complete the proof of Theorem \ref{thm:optimal global estimate for inho DG}.

\subsection{Some estimates on polyharmonic maps with Navier boundary condition}
We begin with a simple lemma concerning the a priori estimate of polyharmonic maps with Navier boundary condition.
\begin{lemma}\label{lem: Apriori estimate}
	Let $B_{1}\subset\R^{n}$ be the unit ball and $g\in W^{m-1,q}(B_{1})$ for some $1<q<2$. Then,
	there exists a unique $h\in W^{m-1,q}(B_{1})$ such that
	\begin{equation}
		\begin{cases}
			\De^{m-1}h=0 & \text{in }B_{1},\\
			\De^{i}h=\De^{i}g & \text{on }\pa B_{1},\quad0\le i\le m-2.
		\end{cases}\label{eq: higher order Navier BVP}
	\end{equation}
	Moreover, there exists a constant $C=C(n,m,q)>0$ such that
	\[
	\|h\|_{W^{m-1,q}(B_{1})}\le C\|g\|_{W^{m-1,q}(B_{1})}.
	\]
\end{lemma}
\begin{proof}
	By a standard approximation argument, we may assume that $g\in C^{\wq}(\bar{B}_{1})$.
	Then the existence and uniqueness of $h$ for equation \eqref{eq: higher order Navier BVP}
	can be deduced easily. The point is to deduce the apriori estimate.
	
	We prove it by induction. When $m=2$, this has been proved in \cite{Guo-Xiang-Zheng-2020-Lp}. Assume the lemma holds with $m$ replaced by $m-1$.
	
	Suppose now $h$ solves equation \eqref{eq: higher order Navier BVP}.
	Put $h_{1}=\De h$. Then,
	\[
	\begin{cases}
		\De^{m-2}h_{1}=0 & \text{in }B_{1},\\
		\De^{i}h_{1}=\De^{i}\De g & \text{on }\pa B_{1},\quad0\le i\le m-3.
	\end{cases}
	\]
	Since $\De g\in W^{m-3,q}(B_{1}),$ by induction, we have $h_{1}\in W^{m-3,q}(B_{1})$
	and
	\[
	\|h_{1}\|_{W^{m-3,q}(B_{1})}\le C\|\De g\|_{W^{m-3,q}(B_{1})}.
	\]
	Returning to the equation of $h$, we have
	\[
	\begin{cases}
		\De h=h_{1} & \text{in }B_{1},\\
		h=g & \text{on }\pa B_{1}.
	\end{cases}
	\]
	Thus by the $L^{p}$ regularity theory of elliptic equations, we obtain
	$h\in W^{m-1,q}(B_{1})$ together with the estimate
	\[
	\|h\|_{W^{m-1,q}(B_{1})}\le C\left(\|h_{1}\|_{W^{m-3,q}(B_{1})}+\|g\|_{W^{m-1,q}(B_{1})}\right).
	\]
	Then, combining it with the estimate of $h_{1}$ yields the desired estimate
	in the lemma. This completes the proof.
\end{proof}

Next we use the above lemma to deduce  the following uniform estimate, which will be used in the proof of Theorem \ref{thm:optimal global estimate for inho DG} later. Let $u, A$ be the functions given in Theorem \ref{thm:optimal global estimate for inho DG}. We have proved that $u\in W^{m+1,\frac{2m}{m+1}}_{\loc}$ in the previous section. Hence $A\De u \in W^{m-1,\frac{2m}{m+1}}_{\loc}$.

\begin{lemma}\label{lemma:uniform estimate}
Let $u, A$ be the functions given in Theorem \ref{thm:optimal global estimate for inho DG}. There exists a constant $C=C(m,p)>0$ satisfying the
	following property. For any $B_{R}(z)\subset B_{1}(0)$, let $h\in W^{m-1,q}(B_{R}(z))$,
	$q=\frac{2m}{m+1}$, be the unique solution of the equation
	\begin{equation}
		\begin{cases}
			\De^{m-1}h=0 & \text{in }B_{R}(z),\\
			\De^{i}h=\De^{i}(A\De u) & \text{on }\pa B_{R}(z),\quad0\le i\le m-2.
		\end{cases}\label{eq: equation in an arbitrary ball}
	\end{equation}
	Then, we have
	\begin{equation}\label{eq: uniform estimate for polyharmonic part}
	\|\na^{m-2}h\|_{L^{2p/(2-p)}(B_{R/2}(z))}\le C\left(\|u\|_{W^{m,2}(B_{1})}+\|f\|_{L^{p}(B_{1})}\right).
	\end{equation}
\end{lemma}
\begin{proof}
	Suppose $h$ solves equation (\ref{eq: equation in an arbitrary ball}).
	Put $h_{R}(x)=R^{2}h(z+Rx)$, $A_{R}=A(z+Rx)$, $u_{R}=u(z+Rx)$ for
	$x\in B_{1}=B_{1}(0)$. Then,
	\[
	\begin{cases}
		\De^{m-1}h_{R}=0 & \text{in }B_{1},\\
		\De^{i}h_{R}=\De^{i}(A_{R}\De u_{R}) & \text{on }\pa B_{1},\quad0\le i\le m-2.
	\end{cases}
	\]
	Then by Lemma \ref{lem: Apriori estimate}, $h_{R}\in W^{m-1,q}(B_{1})$
	with
	\[
	\|h_{R}\|_{W^{m-1,q}(B_{1})}\le C_{m}\|A_{R}\De u_{R}\|_{W^{m-1,q}(B_{1})}.
	\]
	Then, using the Sobolev embedding $W^{m-1,q}(B_1)\subset W^{m-2,2}(B_1)$,
	we obtain a constant $C=C(m,p)>0$ such that
	\[
	\|\na^{m-2}h_{R}\|_{L^{2}(B_{1})}\le C\|\na u_{R}\|_{W^{m,q}(B_{1})}\le CR^{\al}\left(\|u\|_{W^{m,2}(B_{1})}+\|f\|_{L^{p}(B_{1})}\right),
	\]
	where $\al=2m(1-1/p)$. Thus, the previous estimate together with \cite[Lemma 6.2]{Gastel-Scheven-2009CAG} implies that there exists $C=C(m,p)>0$ such that
	\[
	\|\na^{m-2}h_{R}\|_{L^{\bar{p}}(B_{1/2})}\le C\|\na^{m-2}h_{R}\|_{L^{2}(B_{1})}\le CR^{\al}\left(\|u\|_{W^{m,2}(B_{1})}+\|f\|_{L^{p}(B_{1})}\right),
	\]
	which is equivalent to the estimate in the lemma. The proof is complete.
\end{proof}

Now, we are ready to prove Theorem \ref{thm:optimal global estimate for inho DG}.

\subsection{Proof of Theorem \ref{thm:optimal global estimate for inho DG}}
Set
\[
\bar{q}=\frac{2mp}{2m-(m-1)p},\quad  \bar{p}=\frac{2p}{2-p},\quad1<p<\frac{2m}{2m-1}.
\]
The idea is to establish uniform estimate for
$\|\nabla^{m+1} u\|_{L^q(B_{{1}/{2}})}$ with respect to all $q<\bar{q}$ in terms of $\left(\|f\|_{L^{p}(B_{1})}+\|u\|_{L^{1}(B_{1})}\right)$. The approach of Sharp and Topping \cite{Sharp-Topping-2013-TAMS} does not work. We follow the approach of \cite{Guo-Xiang-Zheng-2020-Lp}.  The proof consists of two steps. In the first step, we shall prove $u\in W^{m+1,\bar{q}}_{\loc}$. By Remark \ref{rmk:on best m+1 regularity}, it suffices to show $u\in W^{m,\bar{p}}_{\loc}$. In the second step, we shall derive the optimal interior estimate \eqref{eq:optimal m+1 order} via the conservation law.

For any
\[
\max\left\{2,{\bar{p}}/{2}\right\}<\ga<\bar{p},
\]
we have, for any $1\le i\le m$,
\[
\na^{i}u\in W^{m-i,\ga}(B^{2m})\subset L^{\ga_{i}^{\ast}}(B^{2m}),\quad \text{with}\quad\frac{1}{\ga_{i}^{\ast}}=\frac{1}{\ga}-\frac{m-i}{2m}.
\]
Let $\ga^{\prime}=\ga/(\ga-1)$. Note that $\bar{p}^{\prime}<\ga^{\prime}<2$.
There holds
\[
W_{0}^{m,\bar{p}^{\prime}}(B_{1/2})\subset L^{p^{\prime}}(B_{1/2}).
\]
Thus $W_{0}^{m,\ga^{\prime}}(B_{1/2})\subset L^{\frac{2\ga^{\prime}}{2-\ga^{\prime}}}(B_{1/2})$
with $\frac{2\ga^{\prime}}{2-\ga^{\prime}}=\left(\frac{1}{\ga^{\prime}}-\frac{1}{2}\right)^{-1}=\left(\frac{1}{2}-\frac{1}{\ga}\right)^{-1}>p^{\prime}$.

There are two cases.

\subsubsection{$m$ is an even integer}

We want to estimate the norm of $\|\na^{m}u\|_{\ga,B_{1/2}}$. By \cite[Propositions A.1 and B.1]{Guo-Xiang-Zheng-2020-Lp}, and a direct iteration, we infer that there exists  $C_{p}>0$ independent of $\ga$, such that
\begin{equation}\label{eq: higher order Lp estimate 2}
	\|\na^{m}u\|_{\ga,B_{1/2}}\le C_{p}\left(\|\De^{m/2}u\|_{\ga,B_{1}}+\|u\|_{1,B_{1}}\right).
\end{equation}

To estimate $\|\De^{m/2}u\|_{\ga,B_{1}}$, we  use the system for $A\De u$. More precisely, by Corollary \ref{coro:equation for A Delta u}, we have
\begin{equation}\label{eq:conservation law by coro 3.3}
\De^{m-1}(A\De u)=\sum_{i=1}^{m-1}\de^{i}\left(\sum_{j=m-i}^{m}\na^{j}A\na^{2m-i-j}u\right)+\de K+Af\qquad \text{ in } B_{1}.
\end{equation}
Split $A\De u$ as $A\De u=v+h$ in $B_{1}$ such
that
\[
\begin{cases}
	\De^{m-1}h=0 & \text{in }B_{1},\\
	\De^{i}h=\De^{i}(A\De u) & \text{on }\pa B_{1},\quad0\le i\le m-2
\end{cases}
\]
and
\[
\begin{cases}
	\De^{m-1}v=\De^{m-1}(A\De u) & \text{in }B_{1},\\
	v=\De^{\frac{m}{2}}v=\cdots=\De^{m-2}v=0 & \text{on }\pa B_{1}.
\end{cases}
\]
Then $v$ satisfies the system
$$\Delta^{m-1}v=\sum_{i=1}^{m-1}\de^{i}\left(\sum_{j=m-i}^{m}\na^{j}A\na^{2m-i-j}u\right)+\de K+Af.$$
\medskip

\textbf{Step 1: a duality argument.}
We first estimate $\|\Delta^{\frac{m-2}{2}}v\|_{\gamma,B_1}$.
By duality, we have
\[
\|\De^{\frac{m-2}{2}}v\|_{\ga,B_{1}}=\sup_{\var\in{\cal A}_{1}}\int_{B_{1}}(\De^{\frac{m-2}{2}}v)\var,
\]
where
\[
{\cal A}_{1}=\left\{ \var\in C_{0}^{\wq}(B_{1},\R^{m}):\|\var\|_{\ga^{\prime}}\le1\right\} .
\]
For any $\var\in{\cal A}_{1}$, let $\Phi\in W^{\frac{m}{2},\ga^{\prime}}(B_{1})$
satisfy
\[
\begin{cases}
	\De^{\frac{m}{2}}\Phi=\var & \text{in }B_{1},\\
	\Phi=\De\Phi=\cdots=\De^{\frac{m}{2}-1}\Phi=0 & \text{on }\pa B_{1}.
\end{cases}
\]
By a similar argument as that of \eqref{eq: higher order Lp estimate 2},  there exists a constant $C_p>0$ (independent of $\ga$) such that
\[
\|\Phi\|_{W^{m,\ga^{\prime}}(B_{1})}\le C_{p}\|\var\|_{\ga^{\prime}}\le C_{p}.
\]
Note that integration by parts gives
\[
\begin{aligned}\int_{B_{1}}\De^{\frac{m-2}{2}}v\var & =\int_{B_{1}}\De^{\frac{m-2}{2}}v\De^{\frac{m}{2}}\Phi\\
	& =\int_{B_{1}}\De^{\frac{m}{2}}v\De^{\frac{m-2}{2}}\Phi+\int_{\pa B_{1}}\left(\De^{\frac{m-2}{2}}v\frac{\pa\De^{\frac{m-2}{2}}\Phi}{\pa\nu}-\frac{\pa\De^{\frac{m-2}{2}}v}{\pa\nu}\De^{\frac{m-2}{2}}\Phi\right)\\
	& =\int_{B_{1}}\De^{\frac{m}{2}}v\De^{\frac{m-2}{2}}\Phi\qquad\text{since }\De^{\frac{m-2}{2}}v=\De^{\frac{m-2}{2}}\Phi=0 \text{ on }\partial B_1\\
	& =\int_{B_{1}}\De^{\frac{m}{2}+1}v\De^{\frac{m}{2}-2}\Phi\qquad\text{since }\De^{\frac{m}{2}}v=\De^{\frac{m}{2}-2}\Phi=0\text{ on }\partial B_1\\
	& =\cdots\\
	& =\int_{B_{1}}\De^{m-1}v\Phi\qquad\text{since }\De^{m-2}v=\Phi=0\text{ on }\partial B_1.
\end{aligned}
\]
Thus, we have
\[
\sup_{\var\in{\cal A}_{1}}\int_{B_{1}}\De^{\frac{m-2}{2}}v\cdot\var\le C\sup_{\Phi\in{\cal A}_{2}}\int_{B_{1}}\De^{\frac{m-2}{2}}v\De^{\frac{m}{2}}\Phi=\sup_{\Phi\in{\cal A}_{2}}\int_{B_{1}}\De^{m-1}v\Phi,
\]
where
\[
{\cal A}_{2}=\left\{ \Phi\in W^{m,\ga^{\prime}}(B_{1}):\|\Phi\|_{W^{m,\ga^{\prime}}(B_{1})}\le 1, \Phi=\De\Phi=\cdots=\De^{\frac{m}{2}-1}\Phi=0 \text{ on }\pa B_{1}\right\} .
\]
\medskip

\textbf{Step 2: estimate $\int_{B_{1}}\De^{m-1}v \Phi$ for  $\Phi\in{\cal A}_{2}$.}
Recall that
\[
\int_{B_{1}}\De^{m-1}v\Phi=\int_{B_{1}}\left\{ \sum_{i=1}^{m-1}\de^{i}\left(\sum_{j=m-i}^{m}\na^{j}A\na^{2m-i-j}u\right)+\de K+Af\right\} \Phi.
\]

\textbf{Part 1. }For any $1\le i\le m-1,m-i\le j\le m$, we have
\[
\begin{aligned}
	\int_{B_{1}}\na^{j}A\na^{2m-i-j}u\na^{i}\Phi&\le\|\na^{j}A\|_{\frac{2m}{j}}\left\Vert \na^{2m-i-j}u\right\Vert _{\ga_{2m-i-j}^{\ast}}\|\na^{i}\Phi\|_{\ga_{i}^{\prime\ast}}\\
	&\lesssim\ep\left\Vert \na^{2m-i-j}u\right\Vert _{\ga_{2m-i-j}^{\ast},B_{1}}.
\end{aligned}
\]
So
\[
\int_{B_{1}}\sum_{i=1}^{m-1}\de^{i}\left(\sum_{j=m-i}^{m}\na^{j}A\na^{2m-i-j}u\right)\Phi
\lesssim\ep\sum_{i=1}^{m}\|\na^{i}u\|_{\ga_{i}^{\ast},B_{1}}.
\]

\textbf{Part 2. }By H\"older's inequality,
\[
\int_{B_{1}}(Af)\Phi\le\|Af\|_{p}\|\Phi\|_{p^{\prime},B_{1}}.
\]
Since $\ga<\bar{p}$, we have $\ga^{\prime}>\bar{p}^{\prime}$. Note
that $W_{0}^{m,\bar{p}^{\prime}}(B_{1})\subset L^{p^{\prime}}(B_{1})$.
Thus $W_{0}^{m,\ga^{\prime}}(B_{1})\subset L^{\frac{2\ga^{\prime}}{2-\ga^{\prime}}}(B_{1})$
with $\frac{2\ga^{\prime}}{2-\ga^{\prime}}=\left(\frac{1}{\ga^{\prime}}-\frac{1}{2}\right)^{-1}
=\left(\frac{1}{2}-\frac{1}{\ga}\right)^{-1}>p^{\prime}$.
Here we used the assumption $\ga>2$. Thus
\[
\int_{B_{1}}(Af)\Phi\le C_{p}\|f\|_{p,B_{1}}
\]
for some constant $C_{p}$ independent of $\ga$.

\textbf{Part 3.} \textbf{Estimates concerning $K$.}
We write the first term of $K$ as
\[
K_{1}=\sum_{k=0}^{m-1}\sum_{l=0}^{k-1}(\De^{l}A)\De^{k-l-1}d\langle V_{k},du\rangle
\]
and estimate
\[
\int_{B_{1}}\de K_{1}\Phi=-\int_{B_{1}}K_{1}\cdot\na\Phi
\]
as follows. Repeatedly using the Leibniz rule yields
\[
\De^{l}A\De^{k-l-1}d(V_{k}du)=\sum_{i=1}^{2(k-l)}\De^{l}A\na^{2(k-l)-i}V_{k}\na^{i}u.
\]

\textbf{Case 1:} $2l>m$. Then by \eqref{eq: case 1} there holds
\[
\De^{l}A\na^{2(k-l)-i}V_{k}\na^{i}u  =\sum_{|\al|\le2l-m}\sum_{0\le\be\le\al}\sum_{0\le\eta\le\al-\be}\pa^{\be}\left(A_{\al}\na^{2(k-l)-i+|\eta|}V_{k}\na^{|\al|-|\be|-|\eta|+i}u\right),
\]
where $A_{\al}\in L^{2}(\R^{2m})$. Hence in this case, $K_1$ is a linear combination of terms like
\[
\sum_{|\al|\le2l-m}\sum_{0\le\be\le\al}\sum_{0\le\eta\le\al-\be}\int_{B_{1/2}}A_{\al}\na^{2(k-l)-i+|\eta|}V_{k}\na^{|\al|-|\be|-|\eta|+i}u\na^{|\be|+1}\Phi.
\]
We know  that
\[
\quad\na^{2(k-l)-i-|\eta|}V_{k}\in W^{2l-m+i+1-\eta}\subset L^{\frac{2m}{2m-(2l+i+1)+|\eta|}},\quad\na^{|\al|-|\be|-|\eta|+i}u\in L^{\ga_{|\al|-|\be|-|\eta|+i}^{\ast}}
\]
and so
\[
A_{\al}\na^{2(k-l)-i-|\eta|}V_{k}\na^{|\al|-|\be|-|\eta|+i}u\in L^{q,2}
\]
with
\[
\frac{1}{q}=\frac{1}{2}+\frac{2m-(2l+i+1)+|\eta|}{2m}+\frac{1}{\ga}-\frac{m-(|\al|-|\be|-|\eta|+i)}{2m}=1+\frac{1}{\ga}-\frac{2l-|\al|+|\be|+1}{2m}.
\]
Note that $\frac{2l-|\al|+|\be|+1}{2m}\ge\frac{m+1}{2m}$ for all
$\be,\al$ in the above choice and recall that $\ga>2$. Hence
\[
0<\frac{1}{q}<1+\frac{1}{\ga}-\frac{m+1}{2m}<1.
\]
Also we have
\[
\na^{|\be|+1}\Phi\in W_{0}^{m-|\be|-1,\ga^{\prime}}(B_{1})\subset L^{\frac{2m\ga^{\prime}}{2m-(m-|\be|-1)\ga^{\prime}}}
\]
such that
\[
\frac{1}{q}+\frac{2m-(m-|\be|-1)\ga^{\prime}}{2m\ga^{\prime}}=2-\frac{2l+m-|\al|}{2m}\le1
\]
since $|\al|\le2l-m$. Furthermore, the equality attained only if $|\al|=2l-m$.

Thus
\[
\begin{aligned} & \int_{B_{1/2}}A_{\al}\na^{2(k-l)-i+|\eta|}V_{k}\na^{|\al|-|\be|-|\eta|+i}u\na^{|\be|+1}\Phi\\
	& \le\|A_{\al}\|_{2}\|V_{k}\|\left\Vert \na^{|\al|-|\be|-|\eta|+i}u\right\Vert _{L^{\ga_{|\al|-|\be|-|\eta|+i}^{\ast}}}\left\Vert \na^{|\be|+1}\Phi\right\Vert _{L^{\frac{2m\ga^{\prime}}{2m-(m-|\be|-1)\ga^{\prime}}}}\\
	& \le C_{p}\|A_{\al}\|_{2}\left\Vert \na^{|\al|-|\be|-|\eta|+i}u\right\Vert _{L^{\ga_{|\al|-|\be|-|\eta|+i}^{\ast}}}\\
	& \le C_{p}\ep\left\Vert \na^{|\al|-|\be|-|\eta|+i}u\right\Vert _{L^{\ga_{|\al|-|\be|-|\eta|+i}^{\ast}}}.
\end{aligned}
\]
Summing over all indexes we achieve
\[
\sum_{|\al|\le2l-m}\sum_{0\le\be\le\al}\sum_{0\le\eta\le\al-\be}\int_{B_{1/2}}A_{\al}\na^{2(k-l)-i
+|\eta|}V_{k}\na^{|\al|-|\be|-|\eta|+i}u\na^{|\be|+1}\Phi\lesssim\ep\sum_{i=1}^{m}\|\na^{i}u\|_{\ga_{i}^{\ast},B_{1/2}}.
\]
We can estimate all the remaining terms exactly in the same way as above. Then,
summarizing all the  estimates gives
\begin{equation}\label{eq:estimate for v}
\|\Delta^{\frac{m-2}{2}}v\|_{\gamma,B_1}\lesssim \ep_m\sum_{i=1}^{m}\|\na^{i}u\|_{\ga_{i}^{\ast},B_{1}}+C_{p}\|f\|_{p,B_{1}}.
\end{equation}

\medskip

\textbf{Step 3: Iteration.}
Now, we have
\begin{equation}\label{eq: 6.6}
	\begin{aligned}
		&\qquad \|\nabla^{m-2}(A\Delta u)\|_{L^{\ga}(B_{\frac{1}{2}})}\leq \|\nabla^{m-2}v\|_{L^{\ga}(B_{\frac{1}{2}})}+\|\nabla^{m-2}h\|_{L^\ga(B_{\frac{1}{2}})}\\
		&\leq C\left(\|\Delta^{\frac{m-2}{2}}v\|_{L^{\gamma}(B_1)}+\|u\|_{L^1(B_1)} \right)+\|\nabla^{m-2}h\|_{L^\ga(B_{\frac{1}{2}})}  \\
		&\leq C_p \ep_m\sum_{i=1}^{m}\|\na^{i}u\|_{\ga_{i}^{\ast},B_{1}} + C_p\left(\|f\|_{L^{p}(B_{{1}})}+\|u\|_{L^1(B_1)}+\|\nabla^{m-2}h\|_{L^{\bar{p}}(B_{\frac{1}{2}})}\right).
	\end{aligned}
\end{equation}
Notice the following interpolation inequality from \cite[Chapter 5]{Adams-book} (see also \cite{Sharp-Topping-2013-TAMS}):
\[
\sum_{i=1}^{m-1}\|\na^{i}u\|_{\ga_{i}^{\ast},B_{1}}\le\frac{C}{\ga-1}\left(\|\na^{m}u\|_{\ga,B_{1}}+\|u\|_{1,B_{1}}\right)
\]
for all $\ga>1$, with $C>0$ independent of $\ga$. We obtain from the above interpolation inequality and \eqref{eq: 6.6} that
\begin{equation}\label{eq:key m order estimate}
	\|\na^{m}u\|_{\ga,B_{1/2}}\lesssim \ep_m\|\na^{m}u\|_{\ga,B_{1}}+\|u\|_{1,B_{1}}+\|f\|_{p,B_{1}}+\|\nabla^{m-2}h\|_{\bar{p},B_{\frac{1}{2}}}.
\end{equation}

With \eqref{eq:key m order estimate} at hand, the remaining step is to use a standard scaling technique as that of  \cite[Proof of Lemma 7.2]{Sharp-Topping-2013-TAMS}. Namely, we first use scaling to deduce, for any $B_{R}(z)\subset B_{1}$,
\[
\|\na^{m}u\|_{L^\gamma(B_{\frac{R}{4}}(z))}\le C\ep\|\na^{m}u\|_{L^\ga(B_{R}(z))}+CR^{-2m}(\|f\|_{L^p(B_{R}(z))}
+\|u\|_{L^1(B_{R}(z))}+\|\nabla^{m-2}h\|_{L^{\bar{p}}(B_{\frac{R}{4}})}),
\]
for $\be=2m\bar{p}>0$ (independent of $\ga$). Then, we apply \eqref{eq: uniform estimate for polyharmonic part} to find that
\[
\|\na^{m-2}h\|_{L^{\bar{p}}(B_{\frac{R}{4}}(z))}\le C\left(\|u\|_{W^{m,2}(B_{1})}+\|f\|_{L^{p}(B_{1})}\right)
\]
for all $B_{R}(z)\subset B_{1}$.  Hence
\[
\|\na^{m}u\|_{L^\gamma(B_{\frac{R}{4}}(z))}\le C\ep\|\na^{m}u\|_{L^\ga(B_{R}(z))}+CR^{-2m}\left(\|u\|_{W^{m,2}(B_{1})}+\|f\|_{L^{p}(B_{1})}\right).
\]
Then, we use an iteration lemma of Simon (see e.g. \cite[Lemma A.7]{Sharp-Topping-2013-TAMS}) to derive the following uniform estimate with respect to $\ga$:
\begin{equation}\label{eq: uniform estimate step 1}
	\|\na^{m}u\|_{L^\ga(B_{{1}/{4}})}\le C\left(\|u\|_{W^{m,2}(B_{1})}+\|f\|_{L^{p}(B_{1})}\right)
\end{equation}
with a constant $C=C(p,m)$ independent of $\ga$. Sending $\ga\to\bar{p}$
yields $\na^{m}u\in L^{\bar{p}}(B_{{1}/{4}})$. Consequently, $u\in W^{m+1,\bar{q}}(B_{1/4})$. Thus we  obtain the optimal regularity for $u$.

In the next step, we shall refine estimate \eqref{eq: uniform estimate step 1} to obtain  the following quantitative estimate:
\begin{equation}\label{eq:3}
	\|u\|_{W^{m+1,\bar{q}}(B_{{1}/{2}})}\le C\left(\|f\|_{L^{p}(B_{1})}+\|u\|_{L^{1}(B_{1})}\right).
\end{equation}
For this, we use an interpolation argument.

First of all, note that by the conservation law \eqref{eq:conservation law by coro 3.3},
\begin{equation*}
\De^{m-1}(A\De u)=\sum_{i=1}^{m-1}\de^{i}\left(\sum_{j=m-i}^{m}\na^{j}A\na^{2m-i-j}u\right)+\de K+Af.
\end{equation*}
We may write the right-hand side of the above equation as the form
$\sum_{i=0}^{m-1}\de^if_i$
so as to apply the regularity Lemma \ref{lemma:A2}, where $f_i\in L^{p_i}(B_1)$ are certain vector valued functions with $p_i=\frac{2mp}{2m-ip}$ for all $0\leq i\leq m-1$. Let us consider for instance, the term $\sum_{j=m-i}^{m}\na^{j}A\na^{2m-i-j}u$. Note that $\nabla^jA\in L^{\frac{2m}{j}}$ and $\nabla^{2m-i-j}u\in W^{i+j-m,\bar{p}}\subset L^{\frac{2m\bar{p}}{2m-(i+j-m)\bar{p}}}$. Thus $\na^{j}A\na^{2m-i-j}u\in L^{s}$ with
$$\frac{1}{s}=\frac{j}{2m}+\frac{2m-(i+j-m)\bar{p}}{2m\bar{p}}=\frac{2m-ip}{2mp}=\frac{1}{p_i}.$$
Similarly, we can verify the terms inside $\delta K$ as in the previous calculation. Furthermore, direct computation shows
$$\|f_i\|_{L^{p_i}}\lesssim \|u\|_{W^{m,\bar{p}}}.$$
Combining the above estimate with Lemma \ref{lemma:A2} from the appendix, we conclude
$$\|\nabla^{m+1} u\|_{L^{\bar{q}}(B_{\frac{1}{2}})}\lesssim \|u\|_{W^{m,\bar{p}}(B_1)}+\|u\|_{L^{1}(B_1)}.$$
By standard interpolation and the Sobolev embedding, we have
$$\|u\|_{W^{m,\bar{p}}(B_1)}\leq \epsilon \|\nabla^{m+1} u\|_{L^{\bar{q}}(B_1)}+C_{\epsilon}\|u\|_{L^1(B_1)}.$$
Combining all these estimates, we thus conclude that
\begin{equation}\label{eq:4}
	\|\nabla^{m+1}u\|_{L^{\bar{q}}(B_{\frac{1}{2}})}\lesssim \epsilon \|\nabla^{m+1} u\|_{L^{\bar{q}}(B_1)}+C\left(\|f\|_{L^{p}(B_{1})}+\|u\|_{L^{1}(B_{1})} \right).
\end{equation}
From \eqref{eq:4}, a  scaling argument as that of \eqref{eq: uniform estimate step 1} gives \eqref{eq:3}.

\subsubsection{$m$ is an odd integer}

In the case, $u\in W_{\loc}^{m-1,\eta}$ for any $\eta<\bar{r}=\frac{2mp}{2m-(m+1)p}$ and $m-1$ is an even integer. Furthermore, we have
\[
\|\na^{m-1}u\|_{q,B_{\frac{1}{2}}}\approx\|\De^{\frac{m-3}{2}}(A\Delta u)\|_{q,B_{\frac{1}{2}}}.
\]
In this case, we may repeat exactly what we have done in \textbf{Case I}. The only difference is that instead of first showing that $u\in W^{m,\bar{p}}_{\loc}$, we first show that $u\in W^{m-1,\bar{r}}_{\loc}$ and
$$	\|u\|_{W^{m-1,\bar{r}}(B_{\frac{1}{2}})}\le C\left(\|f\|_{L^{p}(B_{1})}+\|u\|_{W^{m-1,2}(B_{1})}\right).$$
With a similar interpolation argument as in the previous case, we may conclude that
$$\|u\|_{W^{m,\bar{p}}(B_{\frac{1}{2}})}\le C\left(\|f\|_{L^{p}(B_{1})}+\|u\|_{W^{m,2}(B_{1})}\right).$$
%
Then once again an interpolation argument leads to
\begin{equation}\label{eq:before scale and iteration}
	\|\nabla^{m+1}u\|_{\bar{q},B_{\frac{1}{2}}}\leq C_p\ep_m\|\nabla^{m+1}u\|_{\bar{q},B_{{1}}}+C_p\left(\|f\|_{p,B_{1}}+ \|u\|_{1,B_1}\right).
\end{equation}
Finally,  as in the even case dealt above, combining \eqref{eq:before scale and iteration}, a standard scaling argument, together with an iteration lemma of Simon,  gives the desired estimate \eqref{eq:optimal m+1 order}. The proof is thus complete.

\begin{proof}[Proof of Corollary \ref{coro:energy gap}]
	By the previous proof, we know there exists a constant $C=C(p,m)>0$ such that
	$$\sum_{i=1}^m\|\nabla^i u\|_{L^{\bar{q}_i}(B_{\frac{1}{2}})} \leq C(p,m)\|u\|_{L^1(B_1)},$$
	where $\bar{q}_i=\frac{2mp}{2m-(2m-i)p}$, $i=1,\cdots,m$. Using a simple scaling, we then deduce
	\begin{equation}\label{eq:for engergy gap}
		\sum_{i=1}^m\|\nabla^i u\|_{L^{\bar{q}_i}(B_{R})}\leq C(p,m)R^{m(\frac{2-p}{2p}-2)}\|u\|_{W^{m,2}(\R^{2m})}.
	\end{equation}
	Note that $\frac{2-p}{2p}-2<0$. Sending $R\to \infty$ gives $\nabla u=0$ and so $u$ is a constant. Since $u\in L^2(\R^{2m})$, $u\equiv 0$ in $\R^{2m}$.
\end{proof}

\appendix
\section{Higher order $L^{p}$-theory}\label{sec:appendx A}
We write a simple result for higher order $L^p$ regularity theory that's used in the proof of Theorem \ref{thm:optimal global estimate for inho DG}. It is a generalization of the second order equation $-\De u={\rm div}f_1+f_2$.
\begin{lemma} Suppose $n\ge2m$, $1<p<n/m$ and $u\in W^{m,p}(\R^{n})$
solves equation
\[
\De^{m}u=\sum_{i=0}^{m-1}\de^{i}f_{i}\qquad\qquad\text{in }\R^{n}.
\]
where $f_{i}\in L^{p_{i}}(\R^{n})$ with $p_{i}=\frac{np}{n-ip}$
for all $m-1\ge i\ge0$. Then $u\in W^{m+1,\frac{np}{n-(m-1)p}}(\R^{n})$
and
\[
\|u\|_{W^{m+1,\frac{np}{n-(m-1)p}}(\R^{n})}\le C_{m,n,p}\sum_{i=0}^{m-1}\|f_{i}\|_{L^{p_{i}}(\R^{n})}.
\]
\end{lemma}
\begin{proof}
Let $I_{2m}$ be the fundamental solution of $\De^{m}$ in $\R^{n}$.
We have
\[
u=I_{2m}(\sum_{i=0}^{m-1}{\de}^{i}f_{i})\approx\sum_{i=0}^{m-1}I_{2m-i}(f_{i}).
\]
Then $$I_{2m-i}(f_{i})\in W^{2m-i,p_{i}}(\R^{n})\subset W^{m+1,p_{m-1}}(\R^{n})$$
for each $0\le i\le m-1$. Thus
$
u\in W^{m+1,p_{m-1}}(\R^{n}).
$
Moreover, for each $j\le m+1$, there holds
\[
\na^{j}u\approx\sum_{i=0}^{m-1}I_{2m-i-j}(f_{i}).
\]
The Riesz potential theory implies that for each $0\le i\le m-1$,
$I_{2m-i-j}(f_{i})\in L^{q_{j}}(\R^{n})$ with
\[
\frac{1}{q_{j}}=\frac{1}{p_{i}}-\frac{2m-i-j}{n}=\frac{1}{p}-\frac{2m-j}{n}=\frac{n-(2m-j)p}{np},
\]
and
\[
\|\na^{j}u\|_{L^{\frac{np}{n-(2m-j)p}}(\R^{n})}\le C(n,m,p)\sum_{i=0}^{m-1}\|f_{i}\|_{L^{p_{i}}(\R^{n})}.
\]
Therefore,
\[
\sum_{j=1}^{m+1}\|\na^{j}u\|_{L^{\frac{np}{n-(2m-j)p}}(\R^{n})}\le C(n,m,p)\sum_{i=0}^{m-1}\|f_{i}\|_{L^{p_{i}}(\R^{n})}.
\]
This completes the proof.
\end{proof}
Consequently, we can deduce the following local estimate via
a cut-off argument.

\begin{lemma}\label{lemma:A2}
	Suppose $n\ge2m$, $1<p<n/m$ and $u\in W^{m,p}(\R^{n})$
solves equation
\[
\De^{m}u=\sum_{i=0}^{m-1}{\de}^{i}f_{i}\qquad\qquad\text{in }B_{1}.
\]
where $f_{i}\in L^{p_{i}}(B_{1})$ with $p_{i}=\frac{np}{n-ip}$ for
all $m-1\ge i\ge0$. Then
\[
u\in W_{\loc}^{m+1,\frac{np}{n-(m-1)p}}(B_{1})
\]
 and
\[
\|u\|_{W^{m+1,\frac{np}{n-(m-1)p}}(B_{1/2})}\le C_{m,n,p}\left(\sum_{i=0}^{m-1}\|f_{i}\|_{L^{p_{i}}(B_{1})}+\|u\|_{L^{p}(B_{1})}\right).
\]
\end{lemma}

\section{Proof of equation \eqref{eq: 1}}

\begin{proof}
Using induction method. For $m=2$, the claim
holds by a direct computation:
\[
\de\De(Adu)=\de(A\De du+2\na A\cdot\na du+\De Adu)=\de(A\De du)+\de(\na A\na^{2}u+\na^{2}A\na u).
\]

Suppose the claim for $m$ holds. Then, for $m+1$ and $u,A\in W^{m+1,2}$,
we have by induction
\[
\begin{aligned}\de\De^{m}(Adu) & =\De\left(\de(A\De^{m-1}du)+\sum_{i=1}^{m-1}\sum_{j=m-i}^{m}\de^{i}\left(\na^{j}A\na^{2m-i-j}u\right)\right)\\
 & =\de\De\left(A\De^{m-1}du\right)+\sum_{i=1}^{m-1}\sum_{j=m-i}^{m}\de^{i}\De\left(\na^{j}A\na^{2m-i-j}u\right).
\end{aligned}
\]
For the last term, we have
\[
\begin{aligned}\de^{i}\De\left(\na^{j}A\na^{2m-i-j}u\right) & =\de^{i}\left(\na^{j}\De A\na^{2m-i-j}u+\na^{j+1}A\na^{2m-i-j+1}u+\na^{j}A\De\na^{2m-i-j}u\right)\\
 & =\de^{i}\left(\na^{j}\De A\na^{2m-i-j}u+\left(\na^{j}A\De\na^{2m-i-j}u\right)\right)+\de^{i}\left(\na^{j+1}A\na^{2m-i-j+1}u\right)\\
 & =\de^{i+1}\left(\na^{j+1}A\na^{2m-i-j}u\right)-\de^{i}\left(\na^{j+1}\De A\na^{2m-i-j+1}u\right)\\
 & \quad+\de^{i+1}\left(\na^{j}A\na^{2m-i-j+1}u\right)-\de^{i}\left(\na^{j+1}\De A\na^{2m-i-j+1}u\right)\\
 & \quad+\de^{i}\left(\na^{j+1}A\na^{2m-i-j+1}u\right)\\
 & =\de^{i+1}\left(\na^{j+1}A\na^{2m-i-j}u+\na^{j}A\na^{2m-i-j+1}u\right)+\de^{i}\left(\na^{j+1}A\na^{2m-i-j+1}u\right)\\
 & =\de^{i+1}\left(\na^{j+1}A\na^{2(m+1)-(i+1)-(j+1)}u+\na^{j}A\na^{2(m+1)-(i+1)-j}u\right)\\
 & \quad+\de^{i}\left(\na^{j+1}A\na^{2(m+1)-i-(j+1)}u\right).
\end{aligned}
\]
This gives
\[
\sum_{i=1}^{m-1}\sum_{j=m-i}^{m}\de^{i}\De\left(\na^{j}A\na^{2m-i-j}u\right)=\sum_{i=1}^{m+1}\sum_{j=m+1-i}^{m+1}\de^{i}\left(\na^{j}A\na^{2(m+1)-i-j}u\right).
\]
For the term $\de\De\left(A\De^{m-1}du\right)$, it is much easier:
\[
\de\De\left(A\De^{m-1}du\right)=\de\left(A\De^{m}du\right)+\de\left(\na A\De^{m-1}\na^{2}u\right)+\de\left(\De A\De^{m-1}du\right).
\]
The last two terms can be easily disposed by noting that $\De^{m-1}=\de^{m-1}\na^{m-1}$.
Transfer $\de^{m-1}$ from $u$ to $A$ gives
\[
\de\left(\na A\De^{m-1}\na^{2}u\right)+\de\left(\De A\De^{m-1}du\right)=\sum_{i=1}^{m}\de^{i}(\na^{m+1-i}A\na^{m+1}u+\na^{m+2-i}u\na^{m}u).
\]
Hence, in summary, we have
\[
\de\De^{m}(Adu)=\de\left(A\De^{m}du\right)+\sum_{i=1}^{m+1}\sum_{j=m+1-i}^{m+1}\de^{i}\left(\na^{j}A\na^{2(m+1)-i-j}u\right).
\]
The proof is complete.
\end{proof}


\end{document}